\documentclass[11pt]{article}

\usepackage[utf8]{inputenc}
\usepackage[T1]{fontenc}
\usepackage{graphicx}
\usepackage{latexsym}
\usepackage{amsmath,amssymb,amsthm}
\usepackage{mathrsfs}
\usepackage{xcolor}
\usepackage{stmaryrd}
\usepackage{bbm}
\usepackage{url}
\usepackage{longtable}
\usepackage{dsfont}
\usepackage{lmodern}

\usepackage[numbers]{natbib}
\usepackage{color}
\usepackage{lmodern}
\usepackage{float}
\usepackage{enumitem}
\setenumerate{leftmargin=*}
\usepackage{amsopn}
\usepackage{comment}

\usepackage[linktocpage]{hyperref}

\newcommand{\bbn}{\mathbb{N}}

\newcommand{\bbr}{\mathbb{R}}

\newcommand{\bbe}{\mathbb{E}}

\newcommand{\bbp}{\mathbb{P}}

\newcommand{\nto}{{n\to\infty}}

\newcommand{\disfrac}{\displaystyle\frac}

\newcommand{\cals}{{\cal S}}

\newcommand{\calf}{{\cal F}}

\newcommand{\calb}{{\cal B}}

\newcommand{\caln}{{\cal N}}

\newcommand{\al}{{\alpha}}
\newcommand{\la}{{\lambda}}
\newcommand{\La}{{\Lambda}}
\newcommand{\eps}{{\varepsilon}}

\newcommand{\Ga}{{\Gamma}}

\newcommand{\si}{{\sigma}}
\newcommand{\om}{{\omega}}
\newcommand{\Om}{{\Omega}}

\newcommand{\ov}{\overline}

\newcommand{\wt}{\widetilde}

\newcommand{\dd}{\mathrm{d}}

\newcommand{\pd}{\partial}

\newcommand{\Leb}{\mathrm{Leb}}

\newcommand{\lpo}{\preceq}
\newcommand{\gpo}{\succeq}

\DeclareMathOperator*{\esssup}{ess\,sup}
\DeclareMathOperator*{\supp}{supp}

\makeatletter
\newcommand{\opnorm}{\@ifstar\@opnorms\@opnorm}
\newcommand{\@opnorms}[1]{%
	\left|\mkern-1.5mu\left|\mkern-1.5mu\left|
	#1
	\right|\mkern-1.5mu\right|\mkern-1.5mu\right|
}
\newcommand{\@opnorm}[2][]{%
	\mathopen{#1|\mkern-1.5mu#1|\mkern-1.5mu#1|}
	#2
	\mathclose{#1|\mkern-1.5mu#1|\mkern-1.5mu#1|}
}
\makeatother
\newtheoremstyle{neu}
{11pt}      % Space above
{11pt}      % Space below
{}                  % Body font
{}          % Indent amount (empty = no indent, \parindent = para indent)
{\bfseries} % Thm head font
{}          % Punctuation after thm head
{1em}  % Space after thm head: " " = normal interword space;
{\textbf{\thmname{#1}\thmnumber{ #2}\thmnote{ (#3)}.}}          % Thm head spec (can be left empty, meaning 'normal')

\newtheoremstyle{proof}
{11pt}      % Space above
{11pt}      % Space below
{}                  % Body font
{}          % Indent amount (empty = no indent, \parindent = para indent)
{\bfseries} % Thm head font
{}            % Punctuation after thm head
{1em}          % Space after thm head: " " = normal interword space;
{\textbf{\thmname{#1}\thmnote{ #3}.}}          % Thm head spec (can be left empty, meaning 'normal')

\newtheorem{Theorem}{Theorem}[section]
\newtheorem{Corollary}[Theorem]{Corollary}
\newtheorem{Lemma}[Theorem]{Lemma}
\newtheorem{Proposition}[Theorem]{Proposition}
\theoremstyle{neu}
\newtheorem{Definition}[Theorem]{Definition}
\newtheorem{Example}[Theorem]{Example}
\newtheorem{Remark}[Theorem]{Remark}

\theoremstyle{proof}
\newtheorem{Proof}{Proof}

\newcommand{\bthm}{\begin{Theorem}}
	\newcommand{\ethm}{\end{Theorem}}

\newcommand{\bcor}{\begin{Corollary}}
	\newcommand{\ecor}{\end{Corollary}}

\newcommand{\blem}{\begin{Lemma}}
	\newcommand{\elem}{\end{Lemma}}

\newcommand{\bprop}{\begin{Proposition}}
	\newcommand{\eprop}{\end{Proposition}}

\newcommand{\bdf}{\begin{Definition}}
	\newcommand{\edf}{\end{Definition}}

\newcommand{\bex}{\begin{Example}}
	\newcommand{\eex}{\end{Example}}

\newcommand{\brem}{\begin{Remark}}
	\newcommand{\erem}{\end{Remark}}

\newcommand{\bass}{\begin{Assumption}}
	\newcommand{\eass}{\end{Assumption}}

\newcommand{\bpr}{\begin{Proof}}
	\newcommand{\epr}{\end{Proof}}

\newcommand{\benu}{\begin{enumerate}}
	\newcommand{\eenu}{\end{enumerate}}

\newcommand{\bit}{\begin{itemize}}
	\newcommand{\eit}{\end{itemize}}

\newcommand{\lra}{\longrightarrow}
\newcommand{\ind}{\mathbbm{1}}
\newcommand{\bfil}{\boldsymbol{F}}
\newcommand{\ti}{\tilde}
\newcommand{\de}{\delta}
\newcommand{\inj}{\hookrightarrow}

\numberwithin{equation}{section}

\allowdisplaybreaks[3]

\title{Normal approximation of the solution to the stochastic wave equation with Lévy noise}

\date{}

\author{	
	Thomas Delerue\thanks{Chair of Mathematical Statistics, Technical University of Munich, Boltzmannstra\ss e 3, DE-85748 Garching, e-mail: thomas.delerue@tum.de}
}

\begin{document}

	\maketitle

	\begin{abstract}
		For a sequence $\dot{L}^{\eps}$ of Lévy noises with variance $\si^2(\eps)$, we prove the Gaussian approximation of the solution $u^{\eps}$ to the stochastic wave equation driven by $\si^{-1}(\eps) \dot{L}^{\eps}$ and thus extend the result of 
		C. Chong and T. Delerue [\textit{Stoch. Partial Differ. Equ. Anal. Comput.} (2019)] to the class of hyperbolic stochastic PDEs. That is, we find a necessary and sufficient condition in terms of $\si^2(\eps)$ for $u^{\eps}$ to converge in law to the solution to the same equation with Gaussian noise. Furthermore, $u^{\eps}$ is shown to have a space--time version with a càdlàg property determined by the wave kernel, and its derivative $\pd_t u^{\eps}$ a càdlàg version when viewed as a distribution-valued process. These two path properties are essential to our proof of the normal approximation as the limit is characterized by martingale problems that necessitate both random elements. Our results apply to additive as well as to multiplicative noises.		
	\end{abstract}

	\vfill
	
	\noindent
	\begin{tabbing}
		\> secondary: \, 46E35, 60G48
	\end{tabbing}
	
	\vspace{1cm}
	

	\vspace{0.5cm}
	
	\newpage

\section{Introduction}

The wave equation is the prototype of an hyperbolic PDE, widely used e.g. in acoustics and signal processing (\cite{Douglas}). In the literature of stochastic PDEs, the corresponding equation with random perturbation has been extensively studied, especially when the driving noise is Gaussian: See e.g. \cite{Walsh86} for the case of a space--time white noise and Chapter 2 of \cite{Dalang} for a noise that is white in time but spatially correlated. In this paper, we consider the stochastic wave equation 
\begin{equation} \label{WAVE:Levy:intro}
\pd_{tt}  u(t,x) = \pd_{xx} u(t,x) + f(u(t,x)) \dot{L}(t,x), \quad (t,x) \in \bbr^+ \times \bbr,
\end{equation}
where $\dot{L}$ is a Lévy space--time white noise. We investigate the \emph{normal approximation} on $[0,T] \times [0,L]$ of solutions to \eqref{WAVE:Levy:intro} when $\dot{L}$ has a finite second moment and no Gaussian component, that is, when can they be approximated in law on compact domains by the solution to 
\begin{equation} \label{WAVE:Gauss:intro}
\pd_{tt}  u(t,x) = \pd_{xx} u(t,x) + f(u(t,x)) \dot{W}(t,x), \quad (t,x) \in \bbr^+ \times \bbr,
\end{equation}
with a Gaussian space--time white noise $\dot{W}$. The purpose of this work is to show that the \emph{necessary and sufficient} condition for this functional convergence is
\begin{equation} \label{AR:cond}
\lim_{\varepsilon \rightarrow 0} 
\frac{1}{\sigma^2(\varepsilon)}  \int_{\vert z \vert > \kappa \sigma(\varepsilon)} z^2 \, Q^{\varepsilon}( \dd z) = 0
\end{equation}
where $\si^2(\eps)$ is the variance of a homogeneous Lévy noise $\dot{L}^{\eps}$ with Lévy measure $Q^{\varepsilon}$, for all $\kappa > 0$.

Intuitively, to make such an approximation plausible, $\dot{L}^{\eps}$ should be close to $\dot{W}$ in distribution. For Lévy processes having jumps decreasing in size to 0, this was made rigorous in \cite{AR} and more generally in \cite{Cohen07}, where condition \eqref{AR:cond} was first introduced. The passage to an infinite-dimensional setting has been addressed in such generality, to our best knowledge, only in the case of \emph{parabolic} stochastic PDEs like the stochastic heat equation, see \cite{CC:TD} and references therein. We substantially generalize these results to the category of hyperbolic stochastic PDEs. Our second contribution is that we consider throughout equations with \emph{multiplicative noise}. As an application, in the situation of \emph{small jumps approximation} of \cite{AR}, if the impulses of the noise $\dot{L}^{\eps}$ decrease too fast to 0 (such as for a gamma noise), then the corresponding stochastic PDE will not admit a normal approximation, but it will, for instance, if $\dot{L}^{\eps}$ is $\al$-stable for any $\al \in (0,2)$.

The strategy of proof of our main result, Theorem~\ref{main theorem}, is identical to \cite{CC:TD}: We show tightness of the solutions to \eqref{WAVE:Levy:intro} via a generalization of the Aldous criterion \cite{Ivanoff96} and then identify uniquely the limit. Here classical methods relying on the Lévy--Khintchine formula do not apply due to the multiplicative noise, so we resort to \emph{martingale problems} that correspond to the solutions and whose associated martingales \emph{converge} under condition \eqref{AR:cond} to a limit that we link to the solution to \eqref{WAVE:Gauss:intro}. However, their predictable characteristics do not depend on the martingale process itself, which makes other well-established techniques, see e.g. Chapter IX of \cite{JJ}, inapplicable as well. Instead, we prove convergence by hand, so to say, in Skorokhod's representation. 
% To this end, we show that all solution processes considered are supported by a suitable Skorokhod space (and an $L^2$-space as well), see details below.
To this end, we show that all solution processes considered belong to a suitable Skorokhod space (and to an $L^2$-space as well), a fact we will extensively utilize because, in our setting, convergence in Skorokhod topology preserves the martingale property.
%This is another contribution of this paper, that we will extensively utilize because, in our setting, convergence in Skorokhod topology preserves the martingale property.

The \emph{random field solution} $u$ to \eqref{WAVE:Levy:intro} will exhibit a càdlàg property \emph{jointly} in space and time, as we show in Theorem~\ref{cadlag:vers:mildsol}, that is directly linked to the shape of the wave kernel. Now it turns out that to show our normal approximation result, we will need to investigate two different processes simultaneously: $u$ and its time derivative $\pd_t u$, because both appear in the weak formulation of the stochastic wave equation that we consider in this work and, hence, in the aforementioned martingale problems, see Section~\ref{sec:suit:weak:form} for details. This is a substantial difference with \cite{CC:TD} where two different representations of the \emph{same} process needed to be adopted due to the singularities of the heat kernel. We show in Theorem~\ref{theo:cadlag:vers:v:eps} that we can view $\pd_t u$ as a càdlàg process taking values on a space of distributions constructed via Hermite expansions. We also mention that $u$ becomes a \emph{strong martingale} after an appropriate change of coordinate system, a crucial property that we will use to show tightness and the path properties above in place of the factorization method from \cite{DaPrato87, Sanz} (applied in \cite{CC:TD}).

As in \cite{AR, CC:TD, Cohen07}, our motivation comes from numerical simulation: An additional normal approximation of the small jumps of the noise in \eqref{WAVE:Levy:intro} might improve the rate of convergence of numerical schemes, as suggested by the results in \cite{Kohatsu10} for SDEs and in \cite{Chen16} for SPDEs.

This paper is organized as follows. In Section~\ref{sec:prel}, we describe in detail equations \eqref{WAVE:Levy:intro} and \eqref{WAVE:Gauss:intro}. In Section~\ref{sec:funct:set}, we introduce all function spaces needed and show existence of the random elements that will be studied in Section~\ref{sec:main:res}, which contains our main result as well as the main ideas of its proof. The details as well as the proofs for Section~\ref{sec:funct:set} are postponed to Section~\ref{sec:proofs}.

\section{Preliminaries} \label{sec:prel}
Consider on a filtered probability space $(\Om, \calf, \bfil = (\calf_t)_{t \leq T}, \bbp)$ that satisfies the usual conditions, for any  $\eps > 0$, the \emph{stochastic wave equation} on $\bbr^+ \times \bbr$ with vanishing initial conditions:

\begin{equation} \label{WAVE:Levy}
\left[
\begin{array}{ll}
\pd_{tt}  u^{\eps}(t,x) = \pd_{xx} u^{\eps}(t,x) + f(u^{\eps}(t,x))  
\disfrac{\dot{L}^{\eps}(t,x)}{\si(\eps)}, &
(t,x) \in \bbr^+ \times \bbr, \\[1.6mm] \vspace{1.6mm}
u^{\eps}(0,x) = \pd_t u^{\eps}(0, x) = 0, & \textrm{for all} \quad x \in \bbr,
\end{array}	\right.
\end{equation}
where $\dot{L}^{\eps}(t,x)$ is a pure-jump Lévy space--time white noise on $\bbr^+ \times \bbr$ given by 
\begin{equation} \label{Levy}
\begin{split}
L^{\eps}(A) = {}  &  
\int_{\bbr^+ \times \bbr} \int_{\bbr} \ind_{A}(t,x)  z \, 
(\mu^{\eps} - \nu^{\eps})(\dd t, \dd x, \dd z)
\end{split}
\end{equation}
for all bounded Borel sets $A \in \calb_b(\bbr^+ \times \bbr)$. In this representation, $\mu^{\eps}$ is a homogeneous Poisson random measure on $(\bbr^+ \times \bbr) \times \bbr$ relative to the filtration $\bfil$, with intensity measure $\nu^{\eps} = \Leb_{\bbr^+ \times \bbr} \, \otimes \, Q^{\eps}$. Here $Q^{\eps}$ is a Lévy measure on $\bbr$, that is, $Q^{\eps}(\{0\}) = 0$ and $\int_{\bbr} (1 \wedge z^2) \, Q^{\eps}(\dd z) < \infty$, see e.g. Chapter~II in \cite{JJ} for the definition of stochastic integrals with respect to Poisson random measures. Furthermore, we assume that for all $\eps > 0$,
\begin{equation} \label{FinVarLevy}
0<\si^2(\eps) = \int_{\bbr} z^2 \, Q^{\eps}(\dd z) < \infty,
\end{equation}
which is the variance of $L^{\eps}([0,1] \times [0,1])$. The special case
\begin{equation}\label{eq:Qvareps}
Q^{\eps}(A) = \int_{\vert z \vert \leq \eps} \ind_{A}(z) \, Q(\dd z), \quad A \in \calb(\bbr), \quad \eps > 0,
\end{equation}
for a single Poisson random measure $\mu$ having intensity measure $\nu = \Leb_{\bbr^+ \times \bbr} \, \otimes \, Q$, corresponds to the small jump approximation in \cite{AR}.

The function $f \colon \bbr \lra \bbr$ in equation \eqref{WAVE:Levy} will be assumed to be Lipschitz continuous throughout this work.

We are interested in the notion of \emph{mild solution} to \eqref{WAVE:Levy}. It is defined as an $\bfil$-predictable random field $u^{\eps} = \{u^{\eps}(t,x) \mid (t,x) \in \bbr^+ \times \bbr \}$ satisfying for all $(t,x) \in \bbr^+ \times \bbr$,
\begin{equation} \label{Mild:Levy}
\begin{split}
u^{\eps}(t,x) & = 
\int_0^t \int_{\bbr} G_{t-s}(x,y) \frac{f(u^{\eps}(s,y))}{\si(\eps)} \, L^{\eps}(\dd s, \textrm{d}y) \\ & =
\int_{[0,t] \times \bbr} \int_{\bbr} G_{t-s}(x,y)  f(u^{\eps}(s,y)) \frac{z}{\si(\eps)} \,
(\mu^{\eps} - \nu^{\eps})(\dd s, \dd y, \dd z) \quad \mathbb{P}\textrm{-almost surely}. 
\end{split}
\end{equation}
In this equation, $G$ denotes the Green's function of the wave operator
$\pd_{tt} - \pd_{xx}$ and has the following expression:
\begin{equation} \label{Green:Wave}
G_{t-s}(x,y) = G(t,x;s,y) = \frac{1}{2} \ind_{A^+(t,x)}(s,y)
\end{equation}
for any $(t,x,s,y) \in (\bbr^+ \times \bbr)^2$, where 
\begin{equation} \label{back:cone}
A^+(t,x) = \left \{ (s,y) \in \bbr^+ \times \bbr \mid \vert y - x \vert \leq t-s \right \}
\end{equation}
denotes the backward light cone with apex $(t,x)$ restricted to $\bbr^+ \times \bbr$. In particular, $G$ is bounded and not differentiable. By Theorem 3.1 in \cite{CC2}, there exists a unique mild solution $u^{\eps}$ to \eqref{WAVE:Levy} satisfying
\begin{equation} \label{Bound:sol:Levy}
\sup_{(t,x) \in [0,T] \times \bbr} \, \bbe \left[ \big \vert u^{\eps}(t,x) \big \vert^p \right] < \infty
\end{equation}
for all $T > 0$, $0 < p \leq 2$ and $\eps > 0$. Indeed, from \eqref{FinVarLevy} and
\begin{equation} \label{bound:Green}
\sup_{(t,x) \in [0,T] \times \bbr} \, \int_{0}^{t} \int_{\bbr} G(t,x;s,y)^p \, \dd s \, \dd y = (1/2)^p \, T^2 < \infty
\end{equation}
for all $T > 0$ and $p > 0$, Assumption A in \cite{CC2} is easily seen to be satisfied for $p=2$.

We will investigate the \emph{normal approximation} of $u^{\eps}$ and for this, we also consider the solution to the same stochastic PDE as above, but now driven by a Gaussian space--time white noise:
\begin{equation} \label{WAVE:Gauss}
\left[
\begin{array}{ll}
\pd_{tt}  u(t,x) = \pd_{xx} u(t,x) + f(u(t,x)) \dot{W}(t,x), &
(t,x) \in \bbr^+ \times \bbr, \\[1.6mm] \vspace{1.6mm}
u(0,x) = \pd_t u(0, x) = 0, & \textrm{for all} \quad x \in \bbr.
\end{array}	\right.
\end{equation}
The driving noise $\dot{W}$ in \eqref{WAVE:Gauss} is a centered Gaussian random field $\left \{ W(A) \mid A \in \mathcal{B}_b(\bbr^+ \times \bbr) \right \}$ with covariance structure $\mathbb{E}[W(A)W(B)]={\textrm{Leb}}_{\bbr^+ \times \bbr}(A \cap B)$ for bounded Borel sets $A, B \subseteq \bbr^+ \times \bbr$. 

As is well-known (see e.g. Exercise 3.7 of Chapter 3 in \cite{Walsh86}), equation \eqref{WAVE:Gauss} has a (unique) continuous mild solution $u$ that satisfies the corresponding bound in \eqref{Bound:sol:Levy} for all $p > 0$. 
% The bound is automatically gotten in the Picard iteration.

% We will now introduce and describe the function spaces needed in this work to show our normal approximation result. -> link to the introduction

\section{Functional setting} \label{sec:funct:set}

In this work, the letter $C$
%, with subscripts whenever relevant indicating some parameters it depends on,
will always denote a strictly positive constant whose value may change from line to line. Note that $\vert f(x) \vert \leq C \vert x \vert + \vert f(0) \vert $ for all $x \in \bbr$ by the Lipschitz continuity of $f$.

If $\varphi_1, \varphi_2$ are elements of the same $L^2$-space, we will always use the notation $\langle \varphi_1, \varphi_2 \rangle$ for the standard scalar product of that space and $\Vert \cdot \Vert$ for the induced norm. If $\phi$ is an element of a topological vector space and $\phi'$ an element of its topological dual, then $\langle \phi', \phi \rangle$ will always denote the dual pairing of $\phi'$ with $\phi$.

\subsection{Path property of mild solutions} \label{subsec:cadlag:def:and:spaces}
%%The wave kernel (in 1 dimension) has no singularities: It is the indicator function of the \emph{backward} light cone, whose triangular shape fits the definition of $\lpo$.
%This definition fits the triangular shape of the \emph{backward} light cone, whose indicator function is the wave kernel (in 1 dimension). Hence, one can expect the solution to \eqref{WAVE:Levy:intro} not to contain Dirac masses and to exhibit a càdlàg property \emph{jointly} in space and time corresponding to $\lpo$ (loosely speaking, continuity in the \emph{forward} light cone and limits from the flanks; see Definition \ref{def:cadlag} for a precise statement) that, as it turns out, has already been introduced in Section 5 of \cite{CC:VSP}. 

Consider the partial order $\lpo$ on $\bbr^2$:
\begin{equation} \label{part:order}
(\ti{t}, \ti{x}) \lpo (t,x) : \Leftrightarrow \ti{t} \leq t \quad \textrm{and} \quad \vert \ti{x} - x \vert \leq t-\ti{t}
\end{equation}
introduced in Section~5 of \cite{CC:VSP}. We define a space--time càdlàg property corresponding to $\lpo$.
\begin{Definition} \label{def:cadlag}
	A function $\phi \colon M \to \bbr$ with $M \subseteq \bbr^2$ is called \emph{$\lpo$-càdlàg} if for every $(t,x) \in M$,
	\begin{enumerate}
		\item $\lim_{\substack{(\ti{t}, \ti{x}) \to (t, x) \\ (\ti{t}, \ti{x}) \gpo (t, x) }} \phi(\ti{t}, \ti{x}) = \phi(t,x),$
		\item The limits from the \emph{flanks}, that is, 
		\[ \lim_{\substack{(\ti{t}, \ti{x}) \to (t, x) \\ \ti{t} < t, \, \vert \ti{x} - x \vert < t - \ti{t} }} \phi(\ti{t}, \ti{x}), \quad 
		\lim_{\substack{(\ti{t}, \ti{x}) \to (t, x) \\ \ti{x} > x, \,  x - \ti{x} \leq \ti{t} - t < \ti{x} - x }} \phi(\ti{t}, \ti{x}) \quad \textrm{and}  \quad
		\lim_{\substack{(\ti{t}, \ti{x}) \to (t, x) \\ \ti{x} < x, \, \ti{x} - x \leq \ti{t} - t < x - \ti{x} }} \phi(\ti{t}, \ti{x}) \quad \text{all exist.} \]
	\end{enumerate}	
	We further denote the space of all $\lpo$-càdlàg functions on $M$ by $D_{\lpo}(M)$. 
\end{Definition}

We have the following result.
\begin{Theorem} \label{cadlag:vers:mildsol}
	For any $\eps > 0$, let $u^{\eps}$ be a mild solution to the stochastic wave equation \eqref{WAVE:Levy} with noise $\si^{-1}(\eps) \dot{L}^{\eps}$. Then $u^{\eps}$ has a modification $\ov{u}^{\eps}$ in $D_{\lpo}(\bbr^+ \times \bbr)$.
\end{Theorem}

We will investigate the functional convergence of the $\lpo$-càdlàg version $\ov{u}^{\eps}$ of Theorem~\ref{cadlag:vers:mildsol} towards $u$ and to this end, we need a suitable Skorokhod topology for $\lpo$-càdlàg functions. 

Consider the order-preserving change of basis in $\bbr^2$ obtained by rotating the standard basis vectors clockwise by 45 degrees
\begin{equation} \label{basis:change}
H \colon (\bbr^2, \lpo)  \lra (\bbr^2, \leq ), \quad 
\begin{pmatrix} t \\ x \end{pmatrix}  \mapsto 
\frac{1}{\sqrt{2}} \begin{pmatrix}
1 & - 1 \\
1 & 1
\end{pmatrix}
\begin{pmatrix} t \\ x \end{pmatrix}
\end{equation}
as well as a shifting in $\bbr^2$ by $u_0 = (-3/2, 1/2)$ composed with $H$ and then rescaled
\begin{equation} \label{bij:trans:J} 
	J \colon (\bbr^2, \lpo) \lra (\bbr^2, \leq ), \quad u \mapsto \frac{\sqrt{2}}{3} \, H(u - u_0).
\end{equation}
We set $u^* = (3/2, 1/2)$ and use the following notation for closed rectangles with respect to $\lpo$:
\begin{equation} \label{not:2dim:int}
\begin{split}
& [(\ti{t}, \ti{x}), (t,x)]_{\lpo} = \left \{ (s,y) \in \bbr^2 \mid (\ti{t}, \ti{x}) \lpo (s,y) \lpo (t,x) \right \} \quad \textrm{for} \quad (\ti{t}, \ti{x}) \lpo (t,x).
\end{split}
\end{equation}
We then have $[0,1]^2 \subsetneq [u_0, u^*]_{\lpo}$ and $J$ builds a bijection between $[u_0, u^*]_{\lpo}$ and $[0, 1]^2$. This particular choice of the vectors $u_0$ and $u^*$ is for simplicity only, it guarantees that the processes we consider in the proofs of Theorem~\ref{cadlag:vers:mildsol} (and Theorem~\ref{theo:tight:Skor:cadlag:vers:u:eps}) vanish on the axes, a technical requirement of strong martingales often seen in the literature. 

Now let $D([0, 1]^2)$ be the usual Skorokhod space of càdlàg functions on $[0, 1]^2$ with respect to the partial order $\leq$ where $(\ti{t}, \ti{x}) \leq (t,x)$ if and only if $\ti{t} \leq t$ and $\ti{x} \leq x$, see e.g. Section 2 in \cite{Ivanoff80} for a definition. Consider the well-defined bijective transformation
\begin{equation} \label{homeo:trans}
\Phi: D([0, 1]^2) \lra D_{\lpo}([u_0, u^*]_{\lpo}), \quad x \mapsto x \circ J.
\end{equation}
We now draw upon the results of \cite{Straf} on general Skorokhod spaces to obtain the following.

\begin{Lemma} \label{lem:Skor:met}
	 There exists a Skorokhod metric, that will be denoted by $\tau$ throughout this work, that makes $D_{\lpo}([0,1]^2)$ and $D_{\lpo}([u_0, u^*]_{\lpo})$ complete and separable metric spaces and with respect to which the composition
	 \begin{equation} \label{comp:map}
		 D([0, 1]^2) \stackrel{\Phi}{\lra} D_{\lpo}([u_0, u^*]_{\lpo}) \stackrel{\iota}{\hookrightarrow} D_{\lpo} ( [0,1]^2 ),
	 \end{equation}
	 with $\Phi$ as in \eqref{homeo:trans} and $\iota$ the restriction map, is \emph{continuous}. Furthermore, $\lpo$-càdlàg functions are continuous except on at most countably many lines and bounded. If $x_n \stackrel{\tau}{\lra} x$ with $x_n, x \in D_{\lpo} ( [0,1]^2 )$, then $x_n(u) \lra x(u)$ at all continuity points $u$ of $x$.
\end{Lemma}

An immediate consequence of Lemma~\ref{lem:Skor:met} is that \emph{tightness} of probability measures in $D([0,1]^2)$, for which there exist criteria in the literature, implies tightness of the transformed measures (according to \eqref{comp:map}) in $D_{\lpo}( [0,1]^2 )$. This will be of crucial importance for the proof of Theorem~\ref{theo:tight:Skor:cadlag:vers:u:eps}.

% The above lemma shows that some topological problems on $D_{\lpo}( [0,1]^2 )$ may be reduced to problems on the usual Skorokhod space $D( [0, 1]^2 )$, which has already been extensively studied. An important example for this work is \emph{tightness} of probability measures: An immediate consequence of Lemma~\ref{lem:Skor:met} is that tightness in $D( [0,1]^2 )$, for which there exist criteria in the literature, implies tightness in $D_{\lpo}( [0,1]^2 )$. This will be used in Theorem \ref{theo:tight:Skor:cadlag:vers:u:eps}.

A straightforward extension of Lemma~\ref{lem:Skor:met}, in the proof of which a definition of $\tau$ is given, yields a Skorokhod topology on $D_{\lpo}([0,T] \times I)$ with $T > 0$ and $I \subseteq \bbr$ a finite closed interval.

\subsection{Weak formulation} \label{sec:suit:weak:form}

The martingale problem approach mentioned in the introduction relies on a suitable weak formulation of the stochastic wave equation on $\bbr^+ \times \bbr$ that we formally compute from \eqref{WAVE:Levy} in this section. It is inspired by Section~13.1 (together with Definition~9.11) of \cite{Peszat}. 

Let $\phi_1, \phi_2 \in C_c^{\infty}(\bbr)$ and $\eps > 0$. Take the scalar multiplication of both sides of \eqref{WAVE:Levy} with $\phi_2$ and integrate over $[0,t] \times \bbr$. Use the initial condition of 
$\pd_t u^{\eps}$ as well as partial integration twice to obtain
\begin{equation} \label{weak:form:formal:1}
\int_{\bbr} \pd_t u^{\eps}(t,x) \phi_2(x) \, \dd x =
\int_{[0,t] \times \bbr} u^{\eps}(s,x) \phi_2''(x) \, \dd s \, \dd x + 
\int_{0}^{t} \int_{\bbr} \phi_2(x) f(u^{\eps}(s,x)) \disfrac{\dot{L}^{\eps}(s,x)}{\si(\eps)} \, \dd s \, \dd x.
\end{equation}
Now add the equation $\int_{\bbr} u^{\eps}(t,x) \phi_1(x) \, \dd x  = \int_{[0,t] \times \bbr} \pd_t u^{\eps}(s,x) \phi_1(x) \,  \dd s \, \dd x$, which readily follows from the initial condition of $u^{\eps}$, to \eqref{weak:form:formal:1} in order to obtain for all $t \geq 0$,
\begin{equation} \label{weak:form:formal}
\begin{split}
& \int_{\bbr} u^{\eps}(t,x) \phi_1(x) \, \dd x + \int_{\bbr} \pd_t u^{\eps}(t,x) \phi_2(x) \, \dd x \\ & \quad = 
\int_{0}^{t} \left(\int_{\bbr} u^{\eps}(s,x) \phi_2''(x) \, \dd x +  \! \int_{\bbr} \pd_t u^{\eps}(s,x) \phi_1(x) \, \dd x \right) \dd s + \!
\int_{0}^{t} \int_{\bbr} \phi_2(x) f(u^{\eps}(s,x)) \disfrac{\dot{L}^{\eps}(s,x)}{\si(\eps)} \, \dd s \, \dd x.
\end{split}
\end{equation}
It turns out that using equation~\eqref{weak:form:formal:1} alone is enough to prove the necessity of \eqref{AR:cond} for $\ov{u}^{\eps} \stackrel{d}{\lra} u$, but not its sufficiency. For the latter, it is really equation \eqref{weak:form:formal} that will be needed in because it yields an equivalence of weak and mild solution to \eqref{WAVE:Levy} (and analogously for \eqref{WAVE:Gauss}).

Note that in \cite{Walsh86}, page 309, the author develops a different weak formulation for \eqref{WAVE:Gauss}. However, it does not yield an equality of stochastic processes by fixed test function (because the latter must satisfy a condition that depends on the current time point), which is required for martingale problems.

Because $u^{\eps}$ is locally integrable on $\bbr^+ \times \bbr$ by \eqref{Bound:sol:Levy}, it is a distribution on $\bbr^+ \times \bbr$ and so, $\pd_t u^{\eps} = \pd u^{\eps}/ \pd t$ in \eqref{weak:form:formal} will be the time derivative of $u^{\eps}$ in the sense of distributions. The aim of the next section is to find a convenient representation of $\pd u^{\eps}/ \pd t$ by means of a distribution-valued càdlàg process, that we can insert into equation \eqref{weak:form:formal} and thereby use for showing our normal approximation result.

\subsection{Distributional time derivative and path property} \label{subsec:cadlag:def:and:spaces:v:eps}

For simplicity, we write in this paper $\de_{y \pm (t-s)} (\dd x) = \de_{y + (t-s)} (\dd x) + \de_{y - (t-s)} (\dd x)$ and we use this notation for functions as well. Let $\Psi \in C_c^{\infty}(\bbr^+ \times \bbr)$ and fix $(s,y) \in \bbr^+ \times \bbr$. Straightforward computations yield for the Green's function $G$,
\begin{equation*}
\begin{split}
& \int_{\bbr^+ \times \bbr} G(t,x;s,y) \pd_t \Psi(t,x) \, \dd t \, \dd x =
\frac{1}{2} \int_{\bbr} \int_{\bbr^+} 
\ind_{\left \{ \vert y - x \vert \leq t-s \right \}} 
\pd_t \Psi(t,x)  \, \dd t \, \dd x \\ & \qquad =
- \frac{1}{2} \int_{\bbr}
\Psi(s + \vert y - x \vert,x) \, \dd x  = 
- \frac{1}{2} \left(\int_{s}^{\infty} \Psi(t,y + (t-s)) \, \dd t + \int_{s}^{\infty} \Psi(t,y - (t-s)) \, \dd t\right) \\ & \qquad =
- \frac{1}{2} \int_{0}^{\infty} \int_{\bbr} \Psi(t,x) \, \de_{y \pm (t-s)} (\dd x)
\ind_{\{s \leq t\}} \, \dd t.
\end{split}
\end{equation*}
% where we substituted $s + (y-x) = t$ (resp. $s + (x-y) = t$) on $(- \infty, y]$ (resp. $[y, + \infty)$) to obtain the last equality. 
Hence, the distributional time derivative of $G(\cdot, \cdot;s,y)$ on $\bbr^+ \times \bbr$ is a measure on $\bbr^+ \times \bbr$ that we will henceforth denote by $\pd  G/ \pd t \, (\dd t, \dd x; s, y)$ and such that
\begin{equation}  \label{dist:deriv:G}
\frac{\pd  G}{\pd t}(\dd t, \dd x; s, y) = \frac{1}{2} \de_{y \pm (t-s)} (\dd x) \ind_{\{s \leq t\}} \, \dd t =
\frac{\dd  G}{\dd x}(t, \dd x; s, y) \ind_{\{s \leq t\}} \, \dd t
\end{equation}
where $\dd  G/ \dd x \, (t, \dd x; s, y)$ denotes the distributional derivative of $G(t, \cdot;s,y)$ for fixed $t,s,y$, which is readily seen to be equal to $(1/2) \de_{y \pm (t-s)} (\dd x)$ whenever $t \geq s$ and to 0 otherwise.

%We can now proceed to find an expression for $\pd u^{\eps}/ \pd t$. 
Using the expression \eqref{Mild:Levy} of $u^{\eps}$, the stochastic Fubini theorem (see, for example, Theorem 2.6 in \cite{Walsh86}) and \eqref{dist:deriv:G}, we further have
\begin{equation} \label{NR3}
\begin{split}
& \int_{\bbr^+ \times \bbr} u^{\eps}(t,x) \pd_t \Psi(t,x) \, \dd t \, \dd x \\ & \qquad \quad =
\int_{\bbr^+ \times \bbr} 
\left(
\int_{\bbr^+ \times \bbr} G(t,x;s,y) \pd_t \Psi(t,x) \, \dd t \, \dd x
\right)
\frac{f(u^{\eps}(s,y))}{\si(\eps)} \, L^{\eps}(\dd s, \dd y) \\ & \qquad \quad = -
\int_{\bbr^+ \times \bbr} 
\left(
\int_{0}^{\infty} \int_{\bbr} \Psi(t,x) \, 
\frac{\dd  G}{\dd x}(t, \dd x; s, y) \ind_{\{s \leq t\}} \, \dd t
\right)
\frac{f(u^{\eps}(s,y))}{\si(\eps)} \, L^{\eps}(\dd s, \dd y) \\ & \qquad \quad =
- \int_{0}^{\infty} \left(\int_{0}^{t} \int_{\bbr} 
\left(\int_{\bbr} \Psi(t,x) \, \frac{\dd  G}{\dd x}(t, \dd x; s, y)\right)
\frac{f(u^{\eps}(s,y))}{\si(\eps)} \, L^{\eps}(\dd s, \dd y) \right) \dd t \quad \bbp \textrm{-a.s.}
\end{split}
\end{equation}
Recall now that the Schwartz space $\cals(\bbr)$ consists of all $C^{\infty}(\bbr)$-functions with rapid decrease, see e.g. Definition 4.1 in \cite{Demengel}. It has a natural topology induced by the seminorms $\sup_{x \in \bbr} \vert x^j \phi^{(k)}(x) \vert$ with $j, k \in \bbn$ and $\phi \in \cals(\bbr)$. 
%The elements of $\cals'(\bbr)$ are distributions called tempered distributions. 
We define for each $\eps > 0$, an $\cals'(\bbr)$-valued stochastic process
\begin{equation} \label{dist:proc:v}
\begin{split}
v^{\eps}: \bbr^+ & \lra \cals'(\bbr) \\
t & \mapsto \left[ \phi \mapsto 
\int_{0}^{t} \int_{\bbr} 
\left(\int_{\bbr} \phi(x) \, \frac{\dd  G}{\dd x}(t, \dd x; s, y)\right)
\frac{f(u^{\eps}(s,y))}{\si(\eps)} \, L^{\eps}(\dd s, \dd y) \right].
\end{split}
\end{equation}
This stochastic integral is well-defined because by \eqref{Bound:sol:Levy}, the Lipschitz continuity of $f$ and \eqref{dist:deriv:G},
\begin{equation} \label{NR4}
\begin{split}
& \bbe \left[ 
\int_{0}^{t} \int_{\bbr^2} 
\left(\int_{\bbr} \phi(x) \, \frac{\dd  G}{\dd x}(t, \dd x; s, y)\right)^2
\frac{f^2(u^{\varepsilon}(s,y))}{\sigma^2(\varepsilon)}  z^2 \, \nu^{\varepsilon}(\dd s, \textrm{d}y, \textrm{d}z)\right] \\ & \qquad \quad =
\frac{1}{4} \int_0^t \int_{\bbr} \phi^2(y \pm (t-s))
\mathbb{E} \left[ f^2(u^{\varepsilon}(s,y)) \right] \dd s \, \textrm{d}y
\left(\frac{1}{\sigma^2(\varepsilon)} 
\int_{\mathbb{R}} z^2 \, Q^{\varepsilon}(\textrm{d}z)\right) \\ & \qquad \quad \leq C
\int_{0}^{t} \int_{\bbr} \phi^2(y \pm (t-s))  \, \dd s \, \dd y \leq C \int_{\bbr} \phi^2(x) \, \dd x < \infty \quad \textrm{for all} \quad
\phi \in \cals(\bbr)  \quad \textrm{and} \quad t \geq 0.
\end{split}
\end{equation}
Combining the definition of distributions and of $v^{\eps}$, and with \eqref{NR3}, we obtain the following representation for $\pd u^{\eps}/ \pd t$:
\begin{equation} \label{dist:deriv:u:eps}
\left \langle \frac{\pd u^{\eps}}{\pd t}, \Psi \right \rangle = \int_{\bbr^+} \langle v^{\eps}_t, \Psi(t, \cdot) \rangle \, \dd t \quad \textrm{for all} \quad \Psi \in C_c^{\infty}(\bbr^+ \times \bbr).
\end{equation}
Furthermore, using $v^{\eps}$, we can now derive mathematically a weak formulation of \eqref{WAVE:Levy} corresponding to equation \eqref{weak:form:formal}, see Proposition~\ref{prop:mild:implies:weak}.

Actually, $v_t^{\eps}$ is not yet a random distribution: We only have for all Schwartz functions $\phi_1$, $\phi_2$ and scalars $\al_1$, $\al_2$, $\langle v_t^{\eps}, \al_1 \phi_1 + \al_2 \phi_2 \rangle = \al_1 \langle v_t^{\eps}, \phi_1 \rangle + 
\al_2 \langle v_t^{\eps}, \phi_2 \rangle$ $\bbp$-almost surely, so $v_t^{\eps}$ is not a linear functional but rather a random linear functional as defined in \cite{Walsh86} on page 332. 
% A similar problem appears in \eqref{NR3}.
However, we can show that the random field $\{\langle v_t^{\eps}, \phi \rangle \mid \phi \in \cals(\bbr)\}$ has a version with values in $\cals'(\bbr)$. For this, we first recall a few facts on $\cals(\bbr)$. For $q \in \bbn$, let $h_q$ denote the $q$th Hermite function
\begin{equation*}
h_q(x) =  \frac{(-1)^q }{(2^q q! \sqrt{\pi})^{1/2}} e^{x^2/2} \frac{\dd^q}{\dd x^q}e^{-x^2}, \quad x \in \bbr.
\end{equation*} 
As is well-known, $h_q \in \cals(\bbr)$ and a possible orthonormal basis of $L^2(\bbr)$ is given by $\{ h_q \mid q \in \bbn \}$. Define now for each $r \geq 0$, the function space
\begin{equation} \label{def:H:r}
H_r(\bbr) = \left\{ \phi \in L^2(\bbr) \, \, \Big \vert \, \, 
\sum_{q=0}^{\infty} (1 + 2q)^r \langle \phi, h_q \rangle^2 < \infty \right\}.
\end{equation}
Note that this is \emph{not} the fractional Sobolev space on $\bbr$ of order $r$ which is usually defined via the Fourier transform. It is a Hilbert space whose topology is induced by the norm $\Vert \phi \Vert_r = \sqrt{\langle \phi, \phi \rangle_r}$ with the scalar product
\begin{equation} \label{skp:H:r}
\langle \phi, \varphi \rangle_r = 
\sum_{q=0}^{\infty} (1 + 2q)^r \langle \phi, h_q \rangle \langle \varphi, h_q \rangle \quad \textrm{for all} \quad \phi, \varphi \in H_r(\bbr).
\end{equation}
We denote the topological dual of $H_r(\bbr)$ by $H_{-r}(\bbr)$ with dual norm ${\Vert \cdot \Vert}_{-r}$. For each $r \geq 0$ and $q \in \bbn$, consider the continuous and linear functional
\begin{equation} \label{ONB:H-r}
e_{q,-r} \colon H_{r}(\bbr) \lra \bbr, \quad \phi \mapsto (1 +2q)^{-r/2} \langle \phi, h_q \rangle_r.
\end{equation}
By the Riesz representation theorem, the duality $\langle \phi', e_{q,-r} \rangle_{-r} = (1 + 2q)^{-r/2} \langle \phi', h_q \rangle$ holds, the set $\{ e_{q,-r} \mid q \in \bbn \}$ forms an orthonormal basis of $H_{-r}(\bbr)$ and for all $\phi' \in H_{-r}(\bbr)$, we have
\begin{equation} \label{dual:ONB:repr}
\phi' = \sum_{q=0}^{\infty} (1 +2q)^{-r/2} \langle \phi', h_q \rangle e_{q,-r} \quad \textrm{in} \quad H_{-r}(\bbr) \quad \textrm{and} \quad {\Vert \phi' \Vert}^2_{-r} = \sum_{q=0}^{\infty} (1 +2q)^{-r} \langle \phi', h_q \rangle^2.
\end{equation}

By Example 2 in Chapter 4 of \cite{Walsh86}, $\cals(\bbr)$ is then a \emph{nuclear space}, see e.g. pages 330--332 of that chapter for a definition. In particular, $\cals(\bbr) \subseteq H_r(\bbr)$ for all $r \geq 0$ and the injection $(\cals(\bbr), {\Vert \cdot \Vert}_s) \inj (\cals(\bbr), {\Vert \cdot \Vert}_r)$ is a Hilbert-Schmidt operator if $r < s + 1$ (and not  $r < s + 1/2$ as indicated in that example, which is a typo). We then obtain the following regularization.

%Furthermore, for an arbitrary fixed $r > 1$, if we now give the Schwartz space $\cals(\bbr)$ the topology induced by the countable family of norms $\{{\Vert \phi \Vert}_{kr} \mid k \in \bbn \}$, it then becomes a \emph{nuclear space}, see e.g. pages 330--332 of that chapter for a definition. Note that this topology is weaker than the natural topology on $\cals(\bbr)$ introduced earlier. -> Is this the right terminology? \\

%We obtain the chain of inclusions 
%\begin{equation*}
%E = \bigcap_{k \in \bbn} H_{kr}(\bbr) \subseteq H_{r}(\bbr) \subseteq H_0(\bbr) = L^2(\bbr) \subseteq H_{-r}(\bbr) \subseteq \bigcup_{k \in \bbn} H_{-kr}(\bbr) = E'.
%\end{equation*}
%Note that in this setting $E'$ is given the strong topology relative to the norms ${\Vert \phi \Vert}_{kr}$ on $E$, so we do \emph{not} work with the usual topology on $\cals(\bbr)$ nor on $\cals'(\bbr)$ (precise them ? Is $\cals'(\bbr)$ always the same, independently of the topology first chosen in $\cals(\bbr)$?). \\

%+ don't know which additional properties of nuclear spaces I should mention and which no. \\
\begin{Proposition} \label{regula}
	For any $r > 1$, $\eps > 0$ and $t \geq 0$, the random field $\{\langle v_t^{\eps}, \phi \rangle \mid \phi \in \cals(\bbr)\}$ has a version which is in $H_{-r}(\bbr)$ and hence, in $\cals'(\bbr)$ as well.
\end{Proposition}

\bpr
This is a direct application of Theorem 4.1 in \cite{Walsh86}: By the calculation in \eqref{NR4},
\begin{equation*}
\bbe \left[ \Big \vert \langle v_t^{\eps}, \phi \rangle - \langle v_t^{\eps}, \varphi \rangle \Big \vert^2
\right] \leq C \int_{\bbr} \big \vert \phi(x) - \varphi(x) \big \vert^2 \, \dd x \quad \textrm{for all} \quad \phi, \varphi \in \cals(\bbr),
\end{equation*}
so $v_t^{\eps}$ is continuous in probability in the norm ${\Vert \cdot \Vert}_r$ for any $r \geq 0$. \qed
\epr
%+ say that we will from now on keep in mind that $v_t^{\eps}$ is even in $H_{-r}(\bbr)$ for any $r > 1$ thanks to Proposition \ref{regula}

As a consequence, we may and will assume from now on that $v_t^{\eps} \in H_{-r}(\bbr)$ for arbitrary $r > 1$. We then obtain the following path property. 
%with representation \eqref{dist:deriv:u:eps} of $\pd u^{\eps}/ \pd t$ in force. 
\begin{Theorem} \label{theo:cadlag:vers:v:eps}
	For any $r > 2$ and $\eps > 0$, the process $v^{\eps}$ introduced in \eqref{dist:proc:v} has a version $\ov{v}^{\eps}$ in $D(\bbr^+, H_{-r}(\bbr))$, the Skorokhod space of $H_{-r}(\bbr)$-valued càdlàg functions on $\bbr^+$.
\end{Theorem}
In the remainder of this paper, we will work with the càdlàg process $\ov{v}^{\eps}$ obtained in Theorem~\ref{theo:cadlag:vers:v:eps} instead of $\pd u^{\eps}/ \pd t$. Here we point out that even though we will investigate in Section~\ref{sec:main:res} convergence in distribution on \emph{finite} intervals, for technical reasons only (e.g. to avoid tedious calculations related to the boundaries of the interval), we chose in this work a space of distributions on the whole of $\bbr$.

Finally, we follow the same scheme for the continuous mild solution $u$ to \eqref{WAVE:Gauss} and by the proofs of Proposition~\ref{regula} and Theorem~\ref{theo:cadlag:vers:v:eps}, for any $r > 2$, there exists a unique continuous process $\ov{v}$ with values in $H_{-r}(\bbr)$ such that for all $\phi \in \cals(\bbr)$ and $t \geq 0$,
\begin{equation} \label{dist:proc:v:Gauss}
\begin{split}
\langle \ov{v}_t, \phi \rangle & =
\int_{0}^{t} \int_{\bbr} 
\left(\int_{\bbr} \phi(x) \, \frac{\dd  G}{\dd x}(t, \dd x; s, y)\right)
f(u(s,y)) \, W(\dd s, \dd y) \\ & =
\frac{1}{2} \int_{0}^{t} \int_{\bbr} \phi(y \pm (t-s)) f(u(s,y)) \, W(\dd s, \dd y)  \quad  \bbp \textrm{-almost surely}.
\end{split}
\end{equation}

\section{Main result} \label{sec:main:res}

In this section, we fix $T > 0$ as well as $L > 0$. Consider the Cartesian space 
\begin{equation*}
\Om^{\dagger} = \left(D_{\lpo}([0,T] \times [0,L]) \cap L^2([0,T] \times [0,L])\right) \times D([0,T], \bbr).
\end{equation*}
Let $\varrho$ be defined as the sum of the metrics $\tau$ and $d_1$, where $\tau$ is the Skorokhod metric on $D_{\lpo}([0,T] \times [0,L])$, see Lemma~\ref{lem:Skor:met}, and $d_1$ the metric induced by the standard $L^2$-norm on $L^2([0,T] \times [0,L])$. Let also $\tau^{\dagger}$ denote the usual Skorokhod metric on $D([0,T], \bbr)$. We equip $\Om^{\dagger}$ with the product metric
\begin{equation} \label{def:metric:chi:dag}
\chi^{\dagger}((f_1, g_1), (f_2, g_2)) = \varrho(f_1, f_2) + \tau^{\dagger}(g_1, g_2) = \left(\tau(f_1, f_2) + d_1(f_1, f_2)\right) + 
\tau^{\dagger}(g_1, g_2)
\end{equation}
for all $(f_1, g_1), (f_2, g_2) \in \Om^{\dagger}$. The main result of this paper is the following limit theorem.

\begin{Theorem} \label{main theorem}	
	Let $L^{\eps}$ be as in \eqref{Levy} with variance $\si^2(\eps)$ as in \eqref{FinVarLevy} for all $\eps > 0$. Further let $u^{\eps}$ be a mild solution to the stochastic wave equation \eqref{WAVE:Levy} driven by $\si^{-1}(\eps) \dot{L}^{\eps}$, $\ov{u}^{\eps}$ its $\lpo$-càdlàg version given by Theorem~\ref{cadlag:vers:mildsol} and $\ov{v}^{\eps}$ the càdlàg $H_{-r}(\bbr)$-valued process obtained in Theorem~\ref{theo:cadlag:vers:v:eps} for an arbitrary fixed $r > 2$.
	
	In addition, let $u$ be the continuous mild solution to the stochastic wave equation \eqref{WAVE:Gauss} driven by $\dot{W}$ and $\ov{v}$ the continuous $H_{-r}(\bbr)$-valued process satisfying \eqref{dist:proc:v:Gauss}. 
	
	Suppose the Lipschitz function $f$ satisfies $f(0) \neq 0$. We then have
	\begin{equation} \label{main:result}
	(\ov{u}^{\eps}, \langle \ov{v}^{\eps}, \phi \rangle) \stackrel{d}{\lra} (u, \langle \ov{v}, \phi \rangle) \quad \textrm{in} \quad \left( \Om^{\dagger}, \chi^{\dagger} \right) \quad \textrm{as} \quad \eps \to 0 \quad \textrm{for all} \quad \phi \in C_c^{\infty}((0,L))
	\end{equation}	
	if and only if condition \eqref{AR:cond} holds for each $\kappa > 0$.		
\end{Theorem}

\begin{Remark}	
	The weak convergence of $(\langle \ov{u}^{\eps}(t, \cdot), \phi \rangle, \langle \ov{v}^{\eps}_t, \phi \rangle)_{t \le T}$ for all $\phi \in C_c^{\infty}((0,L))$, is needed to show the necessity of \eqref{AR:cond}. That is why the distribution-valued processes $\ov{v}^{\eps}$ and $\ov{v}$ were added to the actual normal approximation of $\ov{u}^{\eps}$ by $u$.	
\end{Remark}

\bpr[of Theorem \ref{main theorem}]
In a first part, we show that \eqref{AR:cond} implies \eqref{main:result}. For any fixed $\phi \in C_c^{\infty}((0,L))$, we will show convergence in distribution of subsequences of $\{(\ov{u}^{\eps}, \langle \ov{v}^{\eps}, \phi \rangle) \mid \eps > 0\}$ toward the limit distribution $(u, \langle \ov{v}, \phi \rangle)$. To this end, we need to consider $\ov{u}^{\eps}$ on the larger domain $[0,T] \times [-T,L+T]$. This is due to the Green's function $G$ (recall \eqref{Green:Wave}): The value of $u(t,x)$ for any $(t,x) \in [0,T] \times [0,L]$ depends on values taken by $u$ on 
$\{(s,x) \in \bbr^+ \times \bbr \mid 0 \leq s \leq t, \, \textrm{dist}(x, [0,L]) \leq t - s \} \subseteq [0,T] \times [-T,L+T]$. This larger domain will be necessary for the proof of Theorem~\ref{theo:weak:implies:mild}. Since we also need to work with $\ov{v}^{\eps}$, in order to prove \eqref{main:result}, we shall work with the second Cartesian space
\begin{equation} \label{def:Om:star}
\begin{split}
\Om^* & = \left(D_{\lpo}\left([0,T] \times [-T, L + T] \right) \cap L^2\left([0,T] \times [-T, L + T] \right)\right) \\ & \qquad \qquad \times 
\left(D([0,T], H_{-r}(\bbr)) \cap L^2([0,T], H_{-r}(\bbr))\right).
\end{split}
\end{equation}
Let $\rho$ be the sum of the usual Skorokhod metric $\tau^*$ on $D([0,T], H_{-r}(\bbr))$ and of the standard $L^2$-metric $d_2$ on $L^2([0,T], H_{-r}(\bbr))$. We then equip $\Om^*$ with the product metric
\begin{equation*}
\chi^*((f_1, g_1), (f_2, g_2)) = \varrho(f_1, f_2) + \rho(g_1, g_2) = \left(\tau(f_1, f_2) + d_1(f_1, f_2) \right) + 
\left( \tau^*(g_1, g_2) + d_2(g_1, g_2) \right)
\end{equation*}
for all $(f_1, g_1), (f_2, g_2) \in \Om^*$.

By Theorem~\ref{theo:tight:Skor:cadlag:vers:u:eps} and Theorem~\ref{theo:tight:L2:cadlag:vers:u:eps}, $\ov{u}^{\eps}$ is tight both in $D_{\lpo}([0,T] \times [-T, L + T])$ and in $L^2([0,T] \times [-T, L + T])$. This readily implies that $\ov{u}^{\eps}$ is also tight in $(D_{\lpo}([0,T] \times [-T, L + T]) \cap L^2([0,T] \times [-T, L + T]), \varrho)$ (it is easy to see that if $K_1$ and $K_2$ are compact sets, one in each function space, then $K_1 \cap K_2$ is compact in the intersection space considered). Analogously, $\ov{v}^{\eps}$ is tight in $(D([0,T], H_{-r}(\bbr)) \cap L^2([0,T], H_{-r}(\bbr)), \rho)$ as a consequence of Theorem~\ref{theo:tight:Skor:cadlag:vers:v:eps} and Corollary~\ref{cor:tight:L2:v:eps}. Since the product of compact spaces is compact, we draw the crucial conclusion that $\{ (\ov{u}^{\eps}, \ov{v}^{\eps}) \mid \eps > 0 \}$ is tight in $(\Om^*, \chi^*)$. Note that no assumptions other than \eqref{Levy} and \eqref{FinVarLevy} on the Lévy noise are needed for this result. 

% We show that the pairs $\{ (\ov{u}^{\eps}, \ov{v}^{\eps}) \mid \eps > 0 \}$ are relatively compact in $(\Om^*, \chi^*)$. Indeed, if $K_1$ and $K_2$ are compact sets, one in each function space, and if, assuming $K_1 \cap K_2 \neq \emptyset$, $(f_n)_{n \in \bbn} \subseteq K_1 \cap K_2$ is a sequence that converges to $f$ in $\tau$ for some $f \in K_1$ and to $f'$ in $d_1$ for some $f' \in K_2$, then by the properties of these two metrics, $f_n(t,x)$ also converges to $f(t,x)$ as well as to $f'(t,x)$ for almost all $(t,x)$. Hence, $f = f'$ almost everywhere on $[0,T] \times [0,L]$, $f_n \stackrel{\varrho}{\lra} f$ and $K_1 \cap K_2$ is (sequentially) compact in the intersection space considered. Analogously, $\ov{v}^{\eps}$ is tight in $(D([0,T], H_{-r}(\bbr)) \cap L^2([0,T], H_{-r}(\bbr)), \rho)$, as a consequence of Theorem \ref{theo:tight:Skor:cadlag:vers:v:eps} and Corollary \ref{cor:tight:L2:v:eps} and since the (deterministic) compactness preservation property above is easily seen to hold in this second intersection space as well. Since the product of compact spaces is compact, we draw the crucial conclusion that $\{ (\ov{u}^{\eps}, \ov{v}^{\eps}) \mid \eps > 0 \}$ is tight in $(\Om^*, \chi^*)$. Note that exactly as in \cite{CC:TD}, no assumptions other than \eqref{Levy} and \eqref{FinVarLevy} on the Lévy noise are needed for this result. 

Subsequently, apply Prokhorov's theorem and let without loss of generality $(\eps_k)_{k \in \bbn}$ be a sequence with $\eps_k \to 0$ such that $(\ov{u}^{\eps_k}, \ov{v}^{\eps_k})_{k \in \bbn}$ converges weakly to some distribution on $(\Omega^*, \chi^*)$ as $k \to \infty$.
%Now $(\Om^*, \chi^*)$ is a complete and separable metric space (note that since  ($\lpo$-)càdlàg functions are bounded, both Skorokhod spaces appearing in $\Om^*$ may be viewed as subspaces of the corresponding $L^2$-spaces, hence separability). 
Then we may further apply Skorokhod's representation theorem, see e.g. Section 1 in \cite{Jakubowski}, and obtain random elements 
\begin{equation} \label{Skor:repr:1}
(w^k, \theta^{k}), (w, \theta) \colon (\ov{\Om}, \ov{\calf}, \ov{\bbp}) \lra (\Om^*, \chi^*),
\end{equation}
defined on a common probability space $(\ov{\Om}, \ov{\calf}, \ov{\bbp})$, and satisfying the following properties:
\begin{equation} \label{Skor:repr:2}
\begin{split}
& (w^k, \theta^{k}) \stackrel{d}{=} 
(\ov{u}^{\eps_k}, \ov{v}^{\eps_k}) \quad \textrm{for all} \quad k \in \bbn \quad \textrm{and} \\ & 
(w^k, \theta^{k})(\ov{\om}) \lra (w, \theta)(\ov{\om}) \quad \textrm{in} \quad (\Om^*, \chi^*) \quad \textrm{as} \quad k \to \infty \quad \textrm{for all} \quad \ov{\om} \in \ov{\Om}.
\end{split}
\end{equation}
We will show in the following that for any $\phi \in C_c^{\infty}((0,L))$,
\begin{equation} \label{ident:limit}
(w, \langle \theta, \phi \rangle) \stackrel{d}{=} (u, \langle \ov{v}, \phi \rangle) \quad \textrm{in} \quad (\Om^{\dagger}, \chi^{\dagger}),
\end{equation}
which together with \eqref{Skor:repr:2} implies
% say what continuous maps are considered? -> No
\begin{equation*}
(\ov{u}^{\eps_k}, \langle \ov{v}^{\eps_k}, \phi \rangle) \stackrel{d}{\lra} (u, \langle \ov{v}, \phi \rangle) \quad \textrm{in} \quad \left( \Om^{\dagger}, \chi^{\dagger} \right) \quad \textrm{as} \quad k \to \infty
\end{equation*}
by the continuous mapping theorem, and altogether, \eqref{main:result}. In this identification step of the distribution of $(w, \langle \theta, \phi \rangle)$, we will refer to the parts of \cite{CC:TD} that are identical.

First, define a filtration $\ov{\bfil} = (\ov{\calf}_t)_{t \leq T}$ on $(\ov{\Om}, \ov{\calf})$:
\begin{equation} \label{def:fil:Om:star}
\ov{\calf}_t = \bigcap_{u \geq t} 
\si \left(w^k(s,x), \theta^{k}_s \mid s \leq u, \, -T \leq x \leq L + T, \, k \in \bbn \right) \vee \caln^{\ov{\bbp}}, \quad 0 \leq t \leq T, 
\end{equation}
where $\caln^{\ov{\bbp}}$ is the set of all $\ov{\bbp}$-null sets of $\ov{\calf}$ (we assume that $\ov{\calf}$ is $\ov{\bbp}$-complete), as well as
\begin{equation} \label{def:char:eta}
\ov{B}_t = \int_{0}^{t} \left(\int_{\bbr} w(s,x) \phi_2''(x) \, \dd x + 
\langle \theta_s, \phi_1 \rangle\right) \dd s \quad \textrm{and} \quad
\ov{C}_t = 
\int_{0}^{t} \int_{\bbr} \phi_2^2(x) f^2(w(s,x)) \, \dd s \, \dd x
\end{equation}
for $\phi_1, \phi_2 \in C_c^{\infty}((-T, L + T))$ and $t \leq T$. 

Assume now for the time being that the pair $(w, \theta)$ satisfies the following martingale problem: The complex-valued càdlàg process
\begin{equation} \label{def:ov:M}
\ov{M}_t = e^{i \xi \left(\langle w(t, \cdot), \phi_1 \rangle + \langle \theta_t, \phi_2 \rangle\right)} - 
\int_{0}^{t} e^{i \xi \left(\langle w(s, \cdot), \phi_1 \rangle + \langle \theta_s, \phi_2 \rangle\right)} \, \ov{A} (\dd s) 
\quad \textrm{with} \quad \ov{A}_t = i \xi \ov{B}_t - \frac{1}{2}  \xi^2 \ov{C}_t, \quad t \leq T,
\end{equation}
is a martingale with respect to $(\overline{\Omega}, \overline{\mathcal{F}}, \overline{\boldsymbol{F}}, \overline{\mathbb{P}})$ for all $\xi \in \mathbb{R}$ and $\phi_1, \phi_2 \in C_c^{\infty}((-T, L + T))$. (Note that because $w$ is $\lpo$-càdlàg, the process $(\langle w(t, \cdot), \phi_1 \rangle)_{t \leq T}$ is càdlàg, and that by a limit argument, $w$, $\theta$ as well as $\ov{M}$ are $\ov{\bfil}$-adapted.) Using \eqref{Skor:repr:2}, we also have
\begin{equation} \label{ess:sup:1}
\esssup_{(t,x) \in [0,T] \times [-T,L+T]} \bbe \left[ \big \vert w(t,x) \big \vert^2 \right] < \infty \quad \textrm{and} \quad  \textrm{for all} \quad x \in \bbr, \quad 
w(0, x) = \theta_0 = 0 \quad \ov{\bbp} \textrm{-a.s.}
\end{equation}
Indeed, the Skorokhod convergence of $w^k$ implies for almost all $(t,x) \in [0,T] \times [-T,L+T]$, $w^k(t,x) \lra w(t,x)$ $\ov{\bbp}$-almost surely and because the projection maps $\pi_{(t,x)}: D_{\lpo}([0,T] \times [-T,L+T]) \lra \bbr$, $f \mapsto f(t,x)$ are measurable, we also have $w^k(t,x) \stackrel{d}{=} \ov{u}^{\eps_k}(t,x)$ for all $(t,x)$. Now the random fields $\{ u^{\eps} \mid \eps >0 \}$  satisfy the uniform bound
\begin{equation} \label{unif:bound:sol:Levy}
\sup_{\eps > 0} \, \sup_{(t,x) \in [0,T] \times \bbr} \,
\mathbb{E}\left[ \vert u^{\varepsilon}(t,x) \vert^2 \right] < \infty,
\end{equation}
which only depends on the Lipschitz function $f$. With \eqref{bound:Green} and \eqref{Bound:sol:Levy}, the proof of \eqref{unif:bound:sol:Levy} goes as Lemma~3.1 in \cite{CC:TD}. (Note that \eqref{unif:bound:sol:Levy} is also crucial for proving the existence of $\ov{u}^{\eps}$ and tightness of $\{ (\ov{u}^{\eps}, \ov{v}^{\eps}) \mid \eps > 0 \}$ in $(\Omega^*, \chi^*)$.) Apply then Fatou's lemma to obtain \eqref{ess:sup:1}.

Assumption \eqref{def:ov:M} enables us, together with \eqref{ess:sup:1}, to show that there exists a Gaussian space--time white noise $\widetilde{W}$ on $[0,T] \times [-T, L + T]$, possibly defined on a complete stochastic basis $(\wt{\Om}, \wt{\calf}, \wt{\bfil}, \wt{\bbp})$ extending $(\ov{\Om}, \ov{\calf}, \ov{\bfil}, \ov{\bbp})$ such that for all $\phi_1, \phi_2 \in C_c^{\infty}((-T,L+T))$,
\begin{equation} \label{weak:form:exact:Gauss}
	\begin{split}
	& \int_{\bbr} w(t,x) \phi_1(x) \, \dd x + \langle \theta_t, \phi_2 \rangle \\ & \qquad \quad = 
	\int_{0}^{t} \left(\int_{\bbr} w(s,x) \phi_2''(x) \, \dd x +  \langle \theta_s, \phi_1 \rangle \right) \dd s + 
	\int_{0}^{t} \int_{\bbr} \phi_2(x) f(w_{-}(s,x)) \, \wt{W}(\dd s, \dd x) \quad \forall t \leq T
	\end{split}
\end{equation}
holds $\wt{\bbp}$-almost surely, where we have set $w_{-}(s,x) = \lim_{r \to s, \, r < s} w(r,x)$. 
% ($w_{-}$ is predictable since $w$ is $\lpo$-càdlàg)
As a consequence, apply Theorem~\ref{theo:weak:implies:mild} to deduce that $w$ is on $[0,T] \times [0,L]$ the continuous mild solution to the stochastic wave equation
\begin{equation} \label{WAVE:Gauss:2}
\left[
\begin{array}{ll}
\pd_{tt}  w(t,x) = \pd_{xx} w(t,x) + f(w(t,x)) \dot{\wt{W}}(t,x), &
(t,x) \in \bbr^+ \times \bbr, \\[1.6mm] \vspace{1.6mm}
w(0,x) = \pd_t w(0, x) = 0, & \textrm{for all} \quad x \in \bbr,
\end{array}	\right.
\end{equation}
of which \eqref{weak:form:exact:Gauss} is the weak formulation (on $[0,T] \times [-T, L + T]$), and that $\theta$ satisfies for all $\phi \in C_c^{\infty}((0,L))$ and $t \leq T$,
\begin{equation} \label{NR18}
\langle \theta_t, \phi \rangle = 
\int_{0}^{t} \int_{\bbr} \left( \int_{\bbr} \phi(x)  
\frac{\dd  G}{\dd x}(t, \dd x; s, y) \right)
f(w(s,y)) \, \wt{W}(\dd s, \dd y) \quad \wt{\bbp} \textrm{-almost surely}.
\end{equation}
Recalling \eqref{dist:proc:v:Gauss}, we then infer that \eqref{ident:limit} holds. 
% Finally, use the proof of Theorem \ref{theo:cadlag:vers:v:eps} to infer that $\langle \theta, \phi \rangle$ is a continuous process for each $\phi \in C_c^{\infty}((0,L))$. 

Now we show how to obtain \eqref{weak:form:exact:Gauss}. Note that $\ov{C}$ in \eqref{def:char:eta} only depends on $\phi_2$, so fix $\phi_1 \in C_c^{\infty}((-T, L + T))$ and apply Theorem~II.2.42 of \cite{JJ} on the process $(\ov{M}_t)_{t \leq T}$ in \eqref{def:ov:M} to first see that $(\ov{B}, \ov{C}, 0)$, with $\ov{B}$ in \eqref{def:char:eta}, are the predictable characteristics of the $\ov{\bfil}$-semimartingale $(\langle w(t, \cdot ), \phi_1 \rangle + \langle \theta_t, \phi_2 \rangle)_{t \leq T}$. Observe that they do not directly depend on the process itself, but on $w$ and $\theta$, which also explains why we are working on the space $(\Om^*, \chi^*)$. Consequently,
\begin{equation*}
M_t(\phi_2) = \int_{\bbr} w(t, x) \phi_1(x) \, \dd x + \langle \theta_t, \phi_2 \rangle - \ov{B}_t, \quad t \leq T,
\end{equation*}
is a continuous square-integrable $\ov{\bfil}$-martingale with quadratic variation process $\ov{C}$ for all $\phi_2 \in C_c^{\infty}((-T, L + T))$. This induces, relative to $(\overline{\Omega}, \overline{\mathcal{F}}, \overline{\boldsymbol{F}}, \overline{\mathbb{P}})$, an orthogonal martingale measure $\left \{ M_t(A), \, t \in [0,T], \, A \in \mathcal{B}([-T, L + T]) \right \}$ (see Chapter 2 in \cite{Walsh86} for a definition) with covariation measure
$Q_M(A \times B \times [s,t]) = \int_{s}^{t} \int_{A \cap B} f^2(w(r,x)) \, \textrm{d}r \, \textrm{d}x$ for all $A, B \in \mathcal{B}([-T, L + T])$. Now use the proof of Theorem~3.13 in \cite{CC:TD} to define, possibly on a complete filtered extension $(\wt{\Om}, \wt{\calf}, \wt{\bfil}, \wt{\bbp})$,
\begin{equation} \label{const:Gauss:noi}
\widetilde{W}_t(\phi_2) \! = \!
\int_{0}^{t}  \int_{\bbr} \ind_{\{ f^2(w_{-}(s,x)) \neq 0 \}} \frac{\phi_2(x)}{f(w_{-}(s,x))}  \, M(\textrm{d}s, \textrm{d}x)  + 
\int_{0}^{t}  \int_{\bbr} \ind_{\{ f^2(w_{-}(s,x)) = 0 \}}  \phi_2(x) \, W'(\textrm{d}s, \textrm{d}x)
\end{equation}
where $W'$ is a Gaussian white noise on $[0,T] \times [-T, L + T]$ independent of $M$, for all $t \leq T$ and $\phi_2 \in C_c^{\infty}((-T, L + T))$, and to further deduce that \eqref{const:Gauss:noi} defines a Gaussian white noise $\wt{W}$ on $[0,T] \times [-T, L + T]$ with respect to $(\wt{\Om}, \wt{\calf}, \wt{\bfil}, \wt{\bbp})$ such that for all $\phi_1, \phi_2 \in C_c^{\infty}((-T,L+T))$, \eqref{weak:form:exact:Gauss} holds $\wt{\bbp}$-almost surely. 

Of course, it remains to show that $\ov{M}$ in \eqref{def:ov:M} is an $\ov{\bfil}$-martingale. To this end, consider first on $\Om$ the càdlàg process 
$(\langle \ov{u}^{\eps}(t, \cdot), \phi_1 \rangle + \langle \ov{v}_t^{\eps}, \phi_2 \rangle)_{t \leq T}$ with $\phi_1, \phi_2 \in C_c^{\infty}((-T,L+T))$ and $\eps > 0$. By Proposition~\ref{prop:mild:implies:weak}, it is indistinguishable from the right-hand side of \eqref{weak:form:exact:Levy}. 
% Recall that $Q^{\eps}$ is the Lévy measure of the noise $L^{\eps}$ and 
Replicate the proof of Theorem~3.8 in \cite{CC:TD} and use \eqref{unif:bound:sol:Levy} to see that for each $\eps > 0$, the pair $(\ov{u}^{\eps}, \ov{v}^{\eps})$ satisfies the following martingale problem: For all $\xi \in \mathbb{R}$ and $\phi_1, \phi_2 \in C_c^{\infty}((-T, L + T))$, the complex-valued process
\begin{equation}  \label{def:M:eps}
\begin{split}
M^{\eps}_t & = e^{i \xi \left(\langle \ov{u}^{\eps}(t, \cdot), \phi_1 \rangle + \langle \ov{v}_t^{\eps}, \phi_2 \rangle\right)} 
- i \xi \int_{0}^{t} e^{i \xi \left( \langle \ov{u}^{\eps}(s, \cdot), \phi_1 \rangle + \langle \ov{v}_s^{\eps}, \phi_2 \rangle \right)}
\left( \langle \ov{u}^{\eps}(s,\cdot), \phi_2'' \rangle +  \langle \ov{v}_s^{\eps}, \phi_1 \rangle \right) \dd s \\ & \quad \,
+ \int_{0}^{t} \int_{\bbr^2} e^{i \xi \left( \langle \ov{u}^{\eps}(s, \cdot), \phi_1 \rangle + \langle \ov{v}_s^{\eps}, \phi_2 \rangle \right)}
\left( e^{i \xi \frac{f(\ov{u}^{\eps}(s,x))}{\si(\eps)} \phi_2(x) z } - 1 - i \xi \frac{f(\ov{u}^{\eps}(s,x))}{\si(\eps)} \phi_2(x) z \right) \dd s \, \dd x \, Q^{\eps}(\dd z)
\end{split}
\end{equation}
is a square-integrable $\bfil$-martingale satisfying $\sup_{\eps > 0} \sup_{t \leq T} \bbe [ \vert M_t^{\eps} \vert^2 ] < \infty$. 
%(Note that because the Lévy noise is homogeneous in space and time, we may switch between $u^{\eps}$ and $\ov{u}^{\eps}$, resp. $v^{\eps}$ and $\ov{v}^{\eps}$, whenever integrating against its intensity measure.)

Define now on $\ov{\Om}$ the $\ov{\bfil}$-adapted process $(\ov{M}_t^k)_{t \leq T}$ in the same way as $M^{\eps}$ in \eqref{def:M:eps}, but with $(\ov{u}^{\eps}, \ov{v}^{\eps})$ and $Q^{\eps}$ replaced by $(w^k, \theta^k)$ of \eqref{Skor:repr:2} and $Q^{\eps_k}$, respectively. Furthermore, because $\ov{M}^k$ has the same distribution as $M^{\eps_k}$ by \eqref{Skor:repr:2}, by standard arguments, $\ov{M}^k$ is a square-integrable $\ov{\bfil}$-martingale satisfying 
$\sup_{k \in \bbn} \sup_{t \leq T} \bbe [ \vert \ov{M}^k_t \vert^2 ] < \infty$ for all $\xi \in \bbr$ and $\phi_1, \phi_2 \in C_c^{\infty}((-T, L + T))$. This is the martingale problem satisfied by the pair $(w^k, \theta^k)$. For any fixed $\xi \in \bbr$ and $\phi_1, \phi_2 \in C_c^{\infty}((-T, L + T))$, we can infer that $\ov{M}$ is an $\ov{\bfil}$-martingale as well, again by standard arguments, if we have:
\begin{equation} \label{conv:ov:M:k:ov:M}
\textrm{for almost all} \quad t \leq T, \quad \ov{M}^k_t \lra \ov{M}_t \quad \textrm{as} \quad k \to \infty \quad \ov{\bbp} \textrm{-almost surely}.
\end{equation}
In order to show \eqref{conv:ov:M:k:ov:M}, which is the final step, first set for each $k \in \bbn$,
\begin{equation} \label{def:char:eta:k}
\begin{split}
& \ov{\nu}^k (A) = 
\int_{0}^{T} \int_{\bbr^2}  
\ind_{A} \left( t, \frac{f(w^k(t,x))}{\si(\eps_k)} \phi_2(x) z \right)  \, \dd t \, \dd x \, Q^{\eps_k}(\dd z),
\\ & \overline{B}_t^{k} = \int_{0}^{t} \left(\langle w^k(s,\cdot), \phi_2'' \rangle +  \langle \theta^k_s, \phi_1 \rangle\right) \textrm{d}s - 
\int_{0}^{t} \int_{\mathbb{R}} x \mathbbm{1}_{\{ \vert x \vert > 1 \}} \, \overline{\nu}^k(\textrm{d} s, \textrm{d} x), \\
& \overline{A}_t^k = i \xi \overline{B}_t^{k} + \int_{0}^{t} \int_{\mathbb{R}} \left( e^{i \xi x} - 1 - i \xi x \mathbbm{1}_{\{ \vert x \vert \leq 1 \}} \right) \,\overline{\nu}^k(\textrm{d}s, \textrm{d}x)
\end{split}
\end{equation}
for all $A \in \mathcal{B}([0,T] \times \mathbb{R})$ and $t \leq T$, so that $\ov{M}^k$ can be written as
\begin{equation} \label{def:ov:M:k}
\begin{split}
\overline{M}_t^k & = e^{i  \xi \left(\langle w^k(t, \cdot), \phi_1 \rangle + \langle \theta_t^k, \phi_2 \rangle\right)}  - 
\int_{0}^{t} e^{i  \xi \left(\langle w^k(s, \cdot), \phi_1 \rangle + \langle \theta_s^k, \phi_2 \rangle\right)} \, \overline{A}^k (\textrm{d}s), \quad  t \leq T.
\end{split}
\end{equation}
Now replicate the proofs of Theorem~3.9, Theorem~3.10 and Lemma~3.11 in \cite{CC:TD} to see that the assumption \eqref{AR:cond} on the Lévy measure $Q^{\eps}$ and \eqref{Skor:repr:2} imply 
\begin{equation} \label{unif:conv}
\sup_{t \leq T} \Bigg \vert 
\int_{0}^{t} e^{i  \xi \left(\langle w^k(s, \cdot), \phi_1 \rangle + \langle \theta_s^k, \phi_2 \rangle\right)} \, \overline{A}^k (\textrm{d}s) - 
\int_{0}^{t} e^{i \xi \left(\langle w(s, \cdot), \phi_1 \rangle + \langle \theta_s, \phi_2 \rangle\right)} \, \ov{A} (\dd s)  \Bigg \vert \lra 0 \quad \textrm{as} \quad k \to \infty
\end{equation}	
pointwise on $\ov{\Om}$. Note that it is the only place in our proof where \eqref{AR:cond} is actually needed. For the proofs of the aforementioned theorems to actually hold here, we need the extra convergence $\int_{0}^{t} \langle \theta^k_s, \phi_1 \rangle \, \dd s \lra \int_{0}^{t} \langle \theta_s, \phi_1 \rangle \, \dd s$. But this readily follows from $\theta^k \lra \theta$ in $L^2([0,T], H_{-r}(\bbr))$ in \eqref{Skor:repr:2}. 
%This explains why the space $L^2([0,T], H_{-r}(\bbr))$ was included in the definition of $\Om^*$ in \eqref{def:Om:star}. Sentence necessary ? -> No, already in introdution.
Recalling the expression of $\ov{M}$ in \eqref{def:ov:M}, resp. of $\ov{M}^k$ in \eqref{def:ov:M:k}, it is now easy to see that \eqref{conv:ov:M:k:ov:M} follows from \eqref{unif:conv}, the Skorokhod convergence of $\theta^k$ and the $L^2$-convergence of $w^k$. \\
%In addition, the Skorokhod convergence of $\theta^k$ implies $\langle \theta_t^k, \phi_2 \rangle \lra \langle \theta_t, \phi_2 \rangle$ $\ov{\bbp}$-almost surely for almost all $t \leq T$, and because $w^k \lra w$ in $L^2([0,T] \times [-T,L+T] )$, we can find a subsequence $(w^{k_l})_{l \in \bbn}$ such that $\langle w^{k_l}(t, \cdot), \phi_1 \rangle \lra \langle w(t, \cdot), \phi_1 \rangle$ for almost all $t \leq T$, pointwise on $\ov{\Om}$. Without loss of generality, we can thus infer
%$\langle w^k(t, \cdot), \phi_1 \rangle + \langle \theta_t^k, \phi_2 \rangle \lra \langle w(t, \cdot), \phi_1 \rangle + \langle \theta_t, \phi_2 \rangle$ $\ov{\bbp}$-almost surely for almost all $t \leq T$, which together with \eqref{unif:conv} and recalling the expression of $\ov{M}$ in \eqref{def:ov:M}, resp. of $\ov{M}^k$ in \eqref{def:ov:M:k}, implies exactly \eqref{conv:ov:M:k:ov:M}. \\

For the second part of the proof, assume \eqref{main:result} and fix a sequence $(\eps_k)_{k \in \bbn}$ converging to 0. We apply again Skorokhod's representation theorem and obtain for each $\phi \in C_c^{\infty}((0,L))$, random elements
\begin{equation*} 
(w^k, \vartheta^{k, \phi}), (w, \vartheta^{\phi}) \colon (\ov{\Om}, \ov{\calf}, \ov{\bbp}) \lra (\Om^{\dagger}, \chi^{\dagger})
\end{equation*}
on a probability space $(\ov{\Om}, \ov{\calf}, \ov{\bbp})$ possibly different from $(\Om, \calf, \bbp)$ but that does \emph{not} depend on $\phi$, satisfying
\begin{equation} \label{Skor:repr:3}
\begin{split}
& (w^k, \vartheta^{k, \phi}) \stackrel{d}{=} 
(\ov{u}^{\eps_k}, \langle \ov{v}^{\eps_k}, \phi \rangle) \quad \textrm{for all} \quad k \in \bbn, \quad 
(w, \vartheta^{\phi}) \stackrel{d}{=} (u, \langle \ov{v}, \phi \rangle) \quad \textrm{and} \\ & 
(w^k, \vartheta^{k, \phi})(\ov{\om}) \lra (w, \vartheta^{\phi})(\ov{\om}) \quad \textrm{in} \quad (\Om^{\dagger}, \chi^{\dagger}) \quad \textrm{as} \quad k \to \infty \quad \textrm{for all} \quad \ov{\om} \in \ov{\Om}.
\end{split}
\end{equation}
We then define a filtration $\ov{\bfil} = (\ov{\calf}_t)_{t \leq T}$ on $(\ov{\Om}, \ov{\calf}, \ov{\bbp})$:
\begin{equation} \label{def:fil:Om:dagger}
\ov{\calf}_t = \bigcap_{u \geq t} 
\si \left(w^k(s,x), \vartheta^{k, \phi}_s \mid s \leq u, \, 0 \leq x \leq L, \, \phi \in C_c^{\infty}((0,L)), \, k \in \bbn \right) \vee \caln^{\ov{\bbp}}, 
\quad 0 \leq t \leq T,
\end{equation}
as well as for arbitrary fixed $\phi\in C_c^{\infty}((0,L))$, the $\ov{\bfil}$-adapted càdlàg processes
\begin{equation} \label{weak:forms:in:ov:Om}
\begin{split}
\ov{X}^k_t & = \vartheta^{k, \phi}_t - \int_{0}^{t} \int_0^L w^k(s,x) \phi''(x) \, \dd x \, \dd s, \\
\ov{X}_t & = \vartheta^{\phi}_t - \int_{0}^{t} \int_0^L w(s,x) \phi''(x) \, \dd x \, \dd s
\end{split}
\end{equation}
for all $k \in \bbn$ and $t \leq T$. Since $\ov{v}$ is continuous, this is also the case for the real-valued process $\vartheta^{\phi}$ by \eqref{Skor:repr:3}, hence $\ov{X}$ is continuous. Now \eqref{Skor:repr:3} readily implies
\begin{equation*}
\ov{X}^k \lra \ov{X} \quad \textrm{in} \quad D([0,T], \bbr) \quad \textrm{as} \quad k \to \infty
\end{equation*}
pointwise on $\ov{\Om}$. Furthermore, by (the proof of) Proposition~\ref{prop:mild:implies:weak}, \eqref{Skor:repr:3} and \eqref{weak:forms:in:ov:Om}, the processes $\ov{X}^k$ and $\ov{X}$ have the same distribution as the square-integrable $\bfil$-martingales
\begin{equation*}
t \mapsto \int_{0}^{t} \int_0^L \phi(x) \disfrac{f(u^{\eps_k}(s,x))}{\si(\eps_k)} \, L^{\eps_k}(\dd s, \dd x) \quad \textrm{and} \quad 
t \mapsto \int_{0}^{t} \int_0^L \phi(x) f(u(s,x)) \, W(\dd s, \dd x),
\end{equation*}
respectively. By standard arguments, we can thus deduce that $\ov{X}^k$ and $\ov{X}$ are $\ov{\bfil}$-martingales. 

Consider the truncation functions
\begin{equation*}
\varrho_h \colon \mathbb{R} \longrightarrow \mathbb{R}, \quad x \mapsto  x \mathbbm{1}_{\{ \vert x \vert \leq h \}}, \quad h > 0,
\end{equation*}
and apply Theorem II.2.21 in~\cite{JJ} to see that the $\ov{\bfil}$-semimartingale characteristics of $\ov{X}^k$ and $\ov{X}$, relative to $\varrho_h$ for a fixed but arbitrary $h > 0$, are given by 
$(\ov{B}^{k,h}, 0, \ov{\nu}^k)$ and $(0, \ov{C}, 0)$, respectively, where $\ov{\nu}^k$ is defined as in \eqref{def:char:eta:k} and $\ov{C}$ as in \eqref{def:char:eta} (with $\phi_2$ replaced by $\phi$), and 
\begin{equation} \label{first semimart charac k}
\overline{B}_t^{k,h} = - \int_{0}^{t} \int_{\mathbb{R}} x \mathbbm{1}_{\{ \vert x \vert > h \}} \, \overline{\nu}^k(\textrm{d} s, \textrm{d} x), \quad t \leq T.
\end{equation}
% + say that the integration limits are automatically changed because of the support of $\phi$? -> No

The remainder of the proof now goes exactly as the proof of Theorem~3.15 in \cite{CC:TD} (where the remaining assumption $f(0) \neq 0$ of the theorem is then needed). \qed
%Note that the only difference is that here the filtration $\ov{\bfil}$ in \eqref{def:fil:Om:dagger} depends on the choice of $\phi$, but this doesn't affect any of the arguments. 
% I changed the filtration -> is it then true?
\epr

\section{Proofs} \label{sec:proofs}

\subsection{Proofs for Section \ref{sec:funct:set} and for tightness}

We begin by showing that each $u^{\eps}$ has a $\lpo$-càdlàg version.

\bpr[of Theorem \ref{cadlag:vers:mildsol}]
Fix $\eps > 0$ for the whole proof. For each $n \in \bbn$, we introduce a truncated Lévy space--time white noise $\dot{L}^{\eps,n}$ on $\bbr^+ \times \bbr$ by setting
\begin{equation}
	\begin{split}
		L^{\eps, n}(A) = {}  &  
		\int_{\bbr^+ \times \bbr} \int_{\bbr} \ind_{A}(t,x) \ind_{\{\vert x \vert \leq n\}} \, 
		z \ind_{\{\vert z \vert > 1/n\}} \, (\mu^{\eps} - \nu^{\eps})(\dd t, \dd x, \dd z)
	\end{split}
\end{equation}
for all $A \in \calb_b(\bbr^+ \times \bbr)$. Now let $u^{\eps,n}$ be a mild solution to the stochastic wave equation \eqref{WAVE:Levy} when $\si^{-1}(\eps) \dot{L}^{\eps}$ is replaced by $\si^{-1}(\eps) \dot{L}^{\eps,n}$. Because $\dot{L}^{\eps,n}$ generates on $[0,t] \times \bbr$ finitely many jumps only, we can write for all $(t,x) \in \bbr^+ \times \bbr$ and $n \in \bbn$,
\begin{equation} \label{sol:trunc:Levy}
	\begin{split}
		u^{\eps, n}(t,x) & = 
		%\int_0^t \int_{\bbr} G_{t-s}(x,y) \frac{f(u^{\eps,n}(s,y))}{\sigma(\eps)} \, L^{\eps,n}(\dd s, \textrm{d}y) \\ & =
		%\int_0^t \int_{\vert y \vert \leq n} \int_{\vert z \vert > 1/n} G_{t-s}(x,y)  f(u^{\eps,n}(s,y)) \frac{z}{\sigma(\eps)} \, \left(\mu^{\eps} (\dd s, \dd y, \dd z) - \nu^{\eps}(\dd s, \dd y, \dd z)\right) \\ & =
		\frac{1}{\sigma(\eps)} \sum_{k=1}^{\infty} 
		G_{t-T_k}(x,X_k)  f(u^{\eps,n}(T_k,X_k)) \ind_{\{\vert X_k \vert \leq n \}} \, Z_k
		\ind_{\{\vert Z_k \vert > 1/n\}} \\ & \quad \,\, - 
		\frac{\int_{\vert z \vert > 1/n} z \, Q^{\eps}(\dd z)}{\si(\eps)} \int_0^t \int_{\vert y \vert \leq n} 
		G_{t-s}(x,y) f(u^{\eps,n}(s,y)) \, \dd s \, \dd y \quad \bbp \textrm{-almost surely}
	\end{split}
\end{equation}
where the $T_k$ indicate the jump times of $\mu^{\eps}$ and $X_k$ (resp. $Z_k$) the space locations (resp. amplitudes) of the jumps of $\mu^{\eps}$. 

The random field on the right-hand side of \eqref{sol:trunc:Levy} is a $\lpo$-càdlàg version of $u^{\eps, n}$. Indeed, through the reformulation of the Green's function
\begin{equation*}
G_{t-s}(x,y) = \frac{1}{2} \ind_{A^-(s,y)}(t,x), \quad (t,x,s,y) \in (\bbr^+ \times \bbr)^2,
\end{equation*}
where 
\begin{equation} \label{forw:cone}
A^-(s,y) = \left \{ (t,x) \in \bbr^+ \times \bbr \mid \vert y - x \vert \leq t-s \right \}
\end{equation}
denotes the forward light cone with apex $(s,y)$, one sees that $(t,x) \mapsto G_{t-T_k}(x,X_k)$ is already $\lpo$-càdlàg and, hence, the finite sum as well as the integral in \eqref{sol:trunc:Levy} are $\bbp$-almost surely $\lpo$-càdlàg. 
% (for the latter, use also , also valid for $u^{\eps,n}$, as well as dominated convergence). 
%the finite sum has an upper index (the number of jumps generated) only depending on $\om \in \Om$ and its summands are all $\lpo$-càdlàg thanks to the $\lpo$-càdlàg property of $(t,x) \mapsto G_{t-T_k}(x,X_k)$ for any fixed $(T_k, X_k)$, as mentioned prior to the statement of Theorem \ref{cadlag:vers:mildsol}. This fact also implies, with \eqref{Bound:sol:Levy} and dominated convergence, that the integral of the second term is almost surely $\lpo$-càdlàg. Hence, each $u^{\eps, n}$ has a $\lpo$-càdlàg version on $\bbr^+ \times \bbr$.

%-> because we are looking for a $\lpo$-càdlàg version of $u^{\eps}$
We will show that $u^{\eps, n}$ converges uniformly on compact sets of $\bbr^+ \times \bbr$ in probability to $u^{\eps}$ as $\nto$. For this, assume first without loss of generality using Theorem~2 in Chapter~3, \S~2 in \cite{Gikhman}, that $u^{\eps}$ and $u^{\eps, n}$ are separable random fields. The first step is to obtain the maximal inequality
\begin{equation} \label{max:in}
	\begin{split}
		& \bbe  \left[\sup_{(s,y) \in [(\ti{t}, \ti{x}), (t,x)]_{\lpo} } \Big \vert u^{\eps}(s,y) - u^{\eps, n}(s,y) \Big \vert^2 \right] \leq
		\bbe  \left[\sup_{(s,y) \in A^+(t,x) } \Big \vert u^{\eps}(s,y) - u^{\eps, n}(s,y) \Big \vert^2 \right]
		 \\ & \qquad \qquad \leq
		\sup_{ (s,y) \in A^+(t,x) }  \bbe  \left[\Big \vert u^{\eps}(s,y) - u^{\eps, n}(s,y) \Big \vert^2 \right] = 
		\bbe \left[ \Big \vert u^{\eps}(t,x) - u^{\eps, n}(t,x) \Big \vert^2 \right]
	\end{split}
\end{equation}
for all $(\ti{t}, \ti{x}) \lpo (t,x)$ in $\bbr^+ \times \bbr$. Choose for simplicity $(\ti{t}, \ti{x}) = 0$ as well as $x = 0$ and fix $t > 0$. Recall the change of coordinates $H$ introduced in \eqref{basis:change} and define $K(u) = H(u - u_0)$ on $\bbr^2$ with $u_0 = (-t,0)$. Then $K$ builds a bijection of $[u_0, u^*]_{\lpo}$ onto $[0,\sqrt{2} t]^2$ with $u^* = (0,t)$. Define also a two-parameter filtration $\bfil^{\eps}$ on $\bbr^2$ with respect to the partial order $\lpo$ by setting
\begin{equation} \label{two:para:fil}
\begin{split}
\calf^{\eps}_{(s, y)} = \bigcap_{(s,y) \lpo (\ti{s}, \ti{y})} \si \left(L^{\eps}(A) \mid A \in \calb( A^+(\ti{s}, \ti{y}) ) \right) \vee \caln^{\bbp} \quad \textrm{for} \quad
(s, y) \in \bbr^+ \times \bbr,
\end{split}
\end{equation}
with $\caln^{\bbp}$ the set of all $\bbp$ null-sets of $\calf$ (and $ \calf^{\eps}_{(s, y)} = \{ \emptyset, \Om \}$ for all $(s,y)$ with $s < 0$). We further define $\wt{u}^{\eps}(v_1,v_2) = u^{\eps} (K^{-1}(v_1, v_2))$ for all $v=(v_1, v_2) \in [0, \sqrt{2} t]^2$ (extending $u^{\eps}$ to 0 whenever $v_1 + v_2 < \sqrt{2} t$) as well as a filtration $\wt{\bfil}^{\eps}$ on $[0,\sqrt{2} t]^2$ with respect to $\leq$ by $\wt{\calf}^{\eps}_{(v_1, v_2)} = \calf^{\eps}_{K^{-1}(v_1, v_2)}$. With the stochastic integration theory of Cairoli and Walsh in \cite{C:W}, we now show that $\wt{u}^{\eps}$ is a two-parameter \emph{strong martingale} with respect to $\wt{\bfil}^{\eps}$, see e.g. page 115 there for a definition. 

Consider on $[0, \sqrt{2} t]^2$ the two-parameter process
\begin{equation*}
\wt{L}^{\eps}(v_1, v_2) = \left\{ \begin{array}{ll}
L^{\eps} \left(A^+(K^{-1}(v_1, v_2)) \right), & \textrm{if} \quad v_1 + v_2 \geq \sqrt{2} t, \\
0, & \textrm{otherwise}.
\end{array} \right.
\end{equation*}
By the properties of the Lévy noise $L^{\eps}$, $\wt{L}^{\eps}$ is a Lévy sheet as well as an $\wt{\bfil}^{\eps}$-strong martingale (the latter follows exactly as in the proof of Lemma~6.2 in \cite{CC:VSP}) and $\wt{\bfil}^{\eps}$ satisfies the commuting condition $F4$ of \cite{C:W}, see pp. 113--114.
%(the difference being that $\wt{L}^{\eps}$ is not trivial on $\{ v_1 + v_2 \geq \sqrt{2} t \}$ only). 
%$\wt{L}^{\eps}$ is actually a Lévy sheet, as is easy to see, and as such has a càdlàg version. 
Choose the filtration $\bfil$ on $[0,t]$ to be 
$\calf_r = \bigcap_{t \geq s \geq r} \si \left(L^{\eps}(A) \mid A \in \calb ( A^+(s,0) ) \right) \vee \caln^{\bbp}$ for all $0 \leq r \leq t$
(note that on $A^+(s,0)$ the mild solution $u^{\eps}$ depends on the values of $L^{\eps}$ on $A^+(s,0)$ only). Then $\wt{u}^{\eps}$ is a valid integrand (see also page 121 of \cite{C:W}) and $\wt{L}^{\eps}$ a valid integrator for Theorem 2.2 in \cite{C:W} to apply, whence
\begin{equation*}
\int_{0}^{v_1} \int_{0}^{v_2} f(\wt{u}^{\eps}(z_1, z_2)) \, \wt{L}^{\eps}(\dd z_1, \dd z_2) =
\int_{\bbr^+ \times \bbr} \ind_{A^+(K^{-1}(v_1, v_2))} f(u^{\eps}(s,y)) \, L^{\eps}(\dd s, \dd y) = \wt{u}^{\eps}(v_1, v_2)
\end{equation*}
is an $\wt{\bfil}^{\eps}$-strong martingale on $[0,\sqrt{2} t]^2$.
%Because on the set $A^+(t,0)$ the mild solution $u^{\eps}$ depends on the values of $L^{\eps}$ on $A^+(t,0)$ only, in the remainder of this proof, we may choose the filtration $\bfil$ on $[0,t]$ to be 
%\begin{equation*}
%\calf_s = \bigcap_{t \geq s \geq r} \si \left(L^{\eps}(A) \mid A \in \calb ( A^+(s,0) ) \right) \vee \caln^{\bbp}, \,\,\, 0 \leq r \leq t.
%\end{equation*}
%Consequently, the push-forward $\wt{u}^{\eps}$ is $\wt{\bfil}^{\eps}$-predictable, see page 121 of \cite{C:W} for the definition of this $\si$-field. Altogether, $\wt{u}^{\eps}$ (resp. $\wt{L}^{\eps}$) is a valid integrand (resp. integrator) and by Theorem 2.2 in \cite{C:W}, the two-parameter integral process
%\begin{equation*}
%\int_{0}^{v_1} \int_{0}^{v_2} f(\wt{u}^{\eps}(z_1, z_2)) \, \wt{L}^{\eps}(\dd z_1, \dd z_2) =
%\int_{\bbr^+ \times \bbr} \ind_{A^+(K^{-1}(v_1, v_2))} f(u^{\eps}(s,y)) \, L^{\eps}(\dd s, \dd y) = \wt{u}^{\eps}(v_1, v_2)
%\end{equation*}
%is an $\wt{\bfil}^{\eps}$-strong martingale on $[0,\sqrt{2} t]^2$.
Analogously, $\wt{u}^{\eps, n} = u^{\eps, n} \circ K^{-1}$ defines an $\wt{\bfil}^{\eps}$-strong martingale on $[0,\sqrt{2} t]^2$ for each $n \in \bbn$.
% (the information of the truncated noise $\dot{L}^{\eps,n}$ is contained in $\bfil^{\eps}$). 
As a consequence, apply Cairoli's strong maximal inequality, see e.g. Corollary~2.3.1 of Chapter~7 in \cite{Khoshnevisan02} (note that $u^{\eps}$ and $u^{\eps,n}$ are also \emph{orthomartingales} by Proposition~1.1 in \cite{Walsh79} and $L^2$-continuous by Theorem~4.7 in \cite{CC2}) to obtain
\begin{equation*}
\begin{split}
& \bbe \left[ \sup_{(v_1, v_2) \in K([0, (t,0)]_{\lpo}) } \Big \vert \wt{u}^{\eps}(v_1, v_2) - \wt{u}^{\eps, n}(v_1, v_2) \Big \vert^2 \right] \leq
\bbe \left[ \sup_{(v_1, v_2) \in [0, \sqrt{2} t]^2 } \Big \vert \wt{u}^{\eps}(v_1, v_2) - \wt{u}^{\eps, n}(v_1, v_2) \Big \vert^2 \right] \\ & \qquad \qquad  \leq 
\sup_{(v_1, v_2) \in [0, \sqrt{2} t]^2} \bbe \left[ \Big \vert \wt{u}^{\eps}(v_1, v_2) - \wt{u}^{\eps, n}(v_1, v_2) \Big \vert^2 \right] = 
\bbe \left[ \Big \vert \wt{u}^{\eps}(\sqrt{2}t,\sqrt{2}t) - \wt{u}^{\eps, n}(\sqrt{2}t,\sqrt{2}t) \Big \vert^2 \right].
\end{split}
\end{equation*}
%+ more precise regarding def/filtrations of orthomartingale and proof of the strong inequality?
By bijectivity, the terms in these inequalities agree exactly with the corresponding ones in \eqref{max:in}.

In a second step, we show that 
\begin{equation} \label{L2:conv}
u^{\eps, n}(t,x) \lra u^{\eps}(t,x) \quad \textrm{in} \quad L^2(\Om, \calf, \bbp) \quad \textrm{as} \quad \nto \quad \textrm{for all} \quad (t,x) \in \bbr^+ \times \bbr.
\end{equation}
Write
\begin{equation*}
	\begin{split}
		u^{\eps}(t,x) - u^{\eps, n}(t,x) & = 
		\int_0^t \int_{\bbr} G_{t-s}(x,y) \frac{f(u^{\eps}(s,y)) - f(u^{\eps, n}(s,y))}{\si(\eps)} \, L^{\eps} (\dd s, \dd y) \\ & \quad + \!
		\int_0^t \int_{\bbr} G_{t-s}(x,y) \frac{f(u^{\eps, n}(s,y))}{\si(\eps)}  (L^{\eps} - L^{\eps,n})(\dd s, \dd y)  =: I_{\eps,n}(t,x) + J_{\eps,n}(t,x).
	\end{split}
\end{equation*}
Fix $T > 0$. Using It\={o}'s isometry and the Lipschitz continuity of $f$, we estimate
\begin{equation} \label{NR2}
	\begin{split}
		\bbe \left[ I_{\eps,n}(t,x)^2 \right] \leq 
		C \int_0^t \int_{\bbr} \ind_{A^+(t,x)}(s,y) \bbe \left[ \Big \vert u^{\eps}(s,y) - u^{\eps,n}(s,y) \Big \vert^2 \right] \dd s \, \dd y
	\end{split}
\end{equation}
as well as
\begin{equation} \label{NR1}
	\begin{split}
		\bbe \left[ J_{\eps,n}(t,x)^2 \right] &  \leq
		C \left(1 + \sup_{(s,y) \in [0,T] \times \bbr} \, \sup_{\eps > 0, n \in \bbn} \, \bbe \left[ \vert u^{\eps,n}(s,y) \vert^2 \right]\right) \\ & \qquad \quad \times \sigma^{-2}(\eps) \int_0^t \int_{\bbr^2} \ind_{A^+(t,x)}(s,y) z^2
		\left(1 - \ind_{\{\vert y \vert \leq n, \vert z \vert > 1/n \} }  \right)^2  \dd s \, \dd y \, Q^{\eps}(\dd z)
	\end{split}
\end{equation}
for all $(t,x) \in [0,T] \times \bbr$ and $n \in \bbn$. Since the uniform bound \eqref{unif:bound:sol:Levy} also applies to all $u^{\eps,n}$ and noting that 
$1 - \ind_{\{\vert y \vert \leq n, \vert z \vert > 1/n \}}  =
\ind_{\{\vert y \vert > n \}} + \ind_{\{\vert y \vert \leq n, \vert z \vert \leq 1/n \}}$ pointwise on $\bbr^2$, the right-hand side of \eqref{NR1} can further be estimated by $C$ times the function
\begin{equation*}
	\begin{split}
		f_{\eps, n}(t,x) & = \int_0^t \int_{\bbr} \ind_{A^+(t,x)}(s,y) \ind_{\{\vert y \vert > n \}} \, \dd s \, \dd y \\ & \quad \,\,  + 
		\sigma^{-2}(\eps) \int_{\bbr} z^2 \ind_{\{\vert z \vert \leq 1/n \}} \, Q^{\eps}(\dd z)  \int_0^t \int_{\bbr}\ind_{A^+(t,x)}(s,y) \, \dd s \, \dd y, \quad 
		(t,x) \in [0,T] \times \bbr,
	\end{split}
\end{equation*}
which together with \eqref{NR2} yields:
\begin{equation} \label{Volterra:ineq}
	\begin{split}
		\bbe \left[ \Big \vert u^{\eps}(t,x) - u^{\eps, n}(t,x) \Big \vert^2 \right] \leq 
		C \int_0^t \int_{\bbr} \ind_{A^+(t,x)}(s,y) \bbe \left[ \Big \vert u^{\eps}(s,y) - u^{\eps,n}(s,y) \Big \vert^2 \right] \dd s \, \dd y + C f_{\eps, n}(t,x)
	\end{split}
\end{equation}
for all $(t,x) \in [0,T] \times \bbr$ and $n \in \bbn$. 

Set $v_{\eps, n}(t,x) = \bbe \left[ \vert u^{\eps}(t,x) - u^{\eps, n}(t,x) \vert^2 \right]$ and hold from now on $C$ in \eqref{Volterra:ineq} fixed. 
%Now $A^+(t,x)$ has Lebesgue measure $t^2$ and $A^+(s,x) \subseteq A^+(t,x)$ if $s \leq t$. 
Let $t_1 > 0$ such that $t_1^2 < 2/C$ and set $t_k = k t_1$ with $k \in \bbn$. We now show by induction that for all $k \in \bbn$, $v_{\eps, n}(t,x) \lra 0$ as $n \to \infty$ for any $(t,x) \in [0,t_k \wedge T] \times \bbr$, which altogether implies \eqref{L2:conv}. First, \eqref{Green:Wave}, \eqref{Volterra:ineq} and dominated convergence yield for $t \leq t_1$,
\begin{equation} \label{ind:beg}
\sup_{(s,y) \lpo (t,x)} v_{\eps, n}(s,y) \leq \frac{C}{1 - C t_1^2/2} \, f_{\eps, n}(t,x) \lra 0 \quad \textrm{as} \quad n \to \infty.
\end{equation}
Next, let $k \geq 2$ and assume $t_k < t \leq t_{k+1} \leq T$. We have
\begin{equation*}
\begin{split}
\int_0^t \int_{\bbr} \ind_{A^+(t,x)}(s,y) v_{\eps, n}(s,y) \, \dd s \, \dd y & \leq
\int_0^{t_k} \int_{\bbr} \ind_{A^+(t,x)}(s,y) v_{\eps, n}(s,y) \, \dd s \, \dd y \\ & \quad \, + 
\sup_{\substack{(\ti{t},\ti{x}) \lpo (t,x) \\ t_k < \ti{t}}} v_{\eps, n}(\ti{t},\ti{x}) 
\int_{t_k}^t \int_{\bbr} \ind_{A^+(t,x)}(s,y) \, \dd s \, \dd y.
\end{split}
\end{equation*}
Combine this inequality with \eqref{Volterra:ineq}, note that $\int_{t_k}^t \int_{\bbr} \ind_{A^+(t,x)}(s,y) \, \dd s \, \dd y = (t- t_k)^2 \leq t_1^2$ and use similar calculations as for \eqref{ind:beg} to conclude that
\begin{equation*}
\begin{split}
\sup_{\substack{(\ti{t},\ti{x}) \lpo (t,x) \\ t_k < \ti{t}}} v_{\eps, n}(\ti{t},\ti{x}) \leq  \frac{C}{1 - C t_1^2/2} \left( \int_0^{t_k} \int_{\bbr} \ind_{A^+(t,x)}(s,y) v_{\eps, n}(s,y) \, \dd s \, \dd y + f_{\eps, n}(t,x) \right) \lra 0
\end{split}
\end{equation*}
as $n \to \infty$ by induction hypothesis and dominated convergence. 

We infer, using \eqref{max:in} and \eqref{L2:conv}, that $u^{\eps, n} - u^{\eps}$ converges uniformly on compacts in probability to 0 as $n \to \infty$ for any $\eps > 0$ and therefore, by standard arguments, the existence of a $\lpo$-càdlàg version $\ov{u}^{\eps}$ of $u^{\eps}$ on $\bbr^+ \times \bbr$.  \qed
%As a consequence, we can finish the proof of the theorem and deduce that each $u^{\eps}$ has a $\lpo$-càdlàg version by the following standard arguments. Take any compact subset $K \subseteq \bbr^+ \times \bbr$. We may assume without loss of generality that $u^{\eps,n} - u^{\eps}$ converges uniformly on $K$ to 0, $\bbp$-almost surely. Consequently, $(u^{\eps,n})_{n \in \bbn}$ is $\bbp$-almost surely a Cauchy sequence with respect to the uniform norm. Since $u^{\eps,n}$ has a $\lpo$-càdlàg version $\ov{u}^{\eps,n}$, as shown in the first part of the proof, and $\lpo$-càdlàg functions are bounded, $\ov{u}^{\eps,n}$ converges uniformly on $K$ to a $\lpo$-càdlàg random field, $\bbp$-almost surely, which is a version of $u^{\eps}$ on $K$. The existence of a $\lpo$-càdlàg version $\ov{u}^{\eps}$ on the whole of $\bbr^+ \times \bbr$ immediately follows. \qed
\epr

We now turn to the Skorokhod topology for $\lpo$-càdlàg functions.

\bpr[of Lemma~\ref{lem:Skor:met}]
We first recall a few facts on the usual Skorokhod topology on $D([0, 1]^2)$ that can all be found in Section 5 of \cite{Straf}. It is induced by the Skorokhod metric
\begin{equation}
\de'(x,y) = \inf_{\la \in \La_s \times \La_s} \left( \sup_{v \in [0, 1]^2} \Big \vert x(v) - y(\la(v)) \Big \vert \vee {\Vert \la \Vert}_s \right), 
\quad x, y \in D([0, 1]^2),
\end{equation}
where $\La_s$ is the set of all homeomorphisms of $[0, 1]$ onto itself which have 0 as a fixed point, $\La_s \times \La_s$ the set of all homeomorphisms $\la$ of the form
\begin{equation*}
\la \colon [0, 1]^2 \lra [0, 1]^2, \quad v=(v_1, v_2) \mapsto (\la_1(v_1), \la_2(v_2))
\end{equation*}
with $\la_1, \la_2 \in \La_s$, and 	
${\Vert \la \Vert}_s =  \sup_{0 \leq p \leq 1}  \left(\max_{i = 1,2} \vert \la_i(p) - p \vert\right)$ for $\la  \in \La_s \times \La_s$. There exists a Skorokhod metric $\de$ that is equivalent to $\de'$ and makes $D([0, 1]^2)$ a complete and separable metric space. 

Now recall \eqref{bij:trans:J}, \eqref{homeo:trans} and give $D_{\lpo}( [0,1]^2)$ the topology induced by the metric
%, thus making $\Phi$ a homeomorphism. 
\begin{equation*}
\tau'(x,y) = \de'(\Phi^{-1}(x), \Phi^{-1}(y)), \quad x, y \in D_{\lpo}([u_0, u^*]_{\lpo}).
\end{equation*}
This is a Skorokhod distance in the sense of \cite{Straf}, see (3.14) of Section~3. Indeed, consider the group of homeomorphisms from $[u_0, u^*]_{\lpo}$ onto itself
$\Theta_s := \{ J^{-1} \circ \la \circ J \mid \la \in \La_s \times \La_s \}$ equipped with the induced norm ${\Vert J^{-1} \circ \la \circ J \Vert}_s := {\Vert \la \Vert}_s$, see Section~3 in \cite{Straf}. Then we can rewrite
\begin{equation*}
\tau'(x,y) = \inf_{\theta \in \Theta_s} \left( \sup_{u \in [u_0, u^*]_{\lpo}} \Big \vert x(u) - 
y(\theta(u)) \Big \vert \vee {\Vert \theta \Vert}_s \right), \quad x, y \in D_{\lpo}([u_0, u^*]_{\lpo}).
\end{equation*}
Defining $\tau$ on $D_{\lpo}([u_0, u^*]_{\lpo})$ analogously to $\tau'$, but with $\de$ instead of $\de'$, yields an equivalent Skorokhod metric to $\tau'$ that makes $D_{\lpo}([u_0, u^*]_{\lpo})$ a complete and separable metric space. 

Definition~\ref{def:cadlag} of $D_{\lpo}([u_0, u^*]_{\lpo})$ coincides exactly with the construction (3.15) in Section~3 of \cite{Straf} of the Skorokhod space on the set $[u_0, u^*]_{\lpo}$ relative to the group $\Theta_s$ (in order to see this, consider all preimages under $J$ of the partitions used in (5.5) and (5.6) of that paper to construct $D([0, 1]^2)$, define the Skorokhod space and use Theorem~5.1 in \cite{Straf}).
% Indeed, note that a collection of finite partitions of $[u_0, u^*]_{\lpo}$ that is invariant under the group action of $\Theta_s$ is needed in that construction: Choose then all preimages under $J$ of the invariant finite partitions used to construct $D([0, 1]^2)$ (see (5.5) and (5.6) in \cite{Straf}), define the Skorokhod space according to (3.15) in \cite{Straf} and use Theorem 5.1 in \cite{Straf}.
%% and consists of uniform limits of sequences of simples functions. As is further shown there, any function of that space may be uniformly approximated by simple functions, which are in our case $\lpo$-càdlàg on $[u_0, u^*]_{\lpo}$ by construction, and therefore is $\lpo$-càdlàg itself. 
% need to be completed?

At last, use the exact same procedure to obtain a Skorokhod topology on $D_{\lpo}\left( [0,1]^2 \right)$ as well as Skorokhod metrics, denoted by the same letters as before. The map $\Phi$ of \eqref{homeo:trans} is now a homeomorphic transformation between 
$D(J([0,1]^2))$ and $D_{\lpo}\left( [0,1]^2 \right)$. For the definition of the Skorokhod metric on $D(J([0,1]^2))$, we now use the subgroup 
$\Ga_s = \{ \la \in \La_s \times \La_s \mid \la(J([0,1]^2)) = J([0,1]^2) \}$ equipped with the norm 
${\Vert \la \Vert}_{s'} = \sup_{(v_1,v_2) \in J([0,1]^2)} \left( \max_{i=1,2} \vert \la_i(v_i) - v_i \vert \right)$. As a consequence, it is easy to see that the restriction map $\iota: D_{\lpo}\left(  [u_0, u^*]_{\lpo} \right) \hookrightarrow D_{\lpo} ( [0,1]^2 )$ is continuous.

The remaining assertions readily follow from Section~3 and~5 of \cite{Straf}. \qed
\epr

Next, we proceed to show that the $\lpo$-càdlàg version $\ov{u}^{\eps}$ is tight. 

\begin{Theorem} \label{theo:tight:Skor:cadlag:vers:u:eps}
	The random fields $\{\ov{u}^{\eps} \mid \eps > 0 \}$ where $\ov{u}^{\eps}$ is the $\lpo$-càdlàg version of $u^{\eps}$ obtained in Theorem~\ref{cadlag:vers:mildsol}, are \emph{tight} in the Skorokhod space $D_{\lpo}\left([0,T] \times I \right)$ for any $T > 0$ and finite closed interval $I \subseteq \bbr$.
\end{Theorem}

\bpr
Without loss of generality, assume that $T=1$ and $I = [0,1]$ and recall the transformation $J$ in \eqref{bij:trans:J}. Set $u^{\eps}(t,x)=\ov{u}^{\eps}(t,x) = 0$ whenever $t < 0$. By Lemma~\ref{lem:Skor:met}, it suffices to show that the random elements $\{ \ov{u}^{\eps} \circ J^{-1} \mid \eps > 0 \}$ are tight in $D([0, 1]^2)$.

By \eqref{unif:bound:sol:Levy}, all random variables $\ov{u}^{\eps} \circ J^{-1}(v_1, v_2)$ are tight. Furthermore, by the same arguments as in the proof of Theorem~\ref{cadlag:vers:mildsol}, the processes $u^{\eps} \circ J^{-1}$ and $\ov{u}^{\eps} \circ J^{-1}$ are strong martingales in $[0, 1]^2$ with respect to the push-forward of filtration \eqref{two:para:fil} through $J$ for each $\eps > 0$. 

We will apply a generalization of Aldous condition for tightness to strong martingales. First of all, fix $\eps > 0$ and note that if $\tau$ is a natural $1$-stopping time for $\ov{u}^{\eps} \circ J^{-1}$ with $\tau \in [0, 1]$, see page 112 in \cite{Ivanoff96} for a definition, then the processes $(\ov{u}^{\eps} \circ J^{-1}(\tau,v))_{0 \leq v \leq 1}$ and 
$(u^{\eps} \circ J^{-1}(\tau,v))_{0 \leq v \leq 1}$ are versions of one another. To see this, approximate $\tau$ from above by a sequence $(\tau_n)_{n \in \bbn}$ of natural $1$-stopping times taking on finitely many values only. Then for all $0 \leq v \leq 1$, $\ov{u}^{\eps} \circ J^{-1}(\tau_n,v) = u^{\eps} \circ J^{-1}(\tau_n,v)$ 
$\bbp$-almost surely. Now since $\ov{u}^{\eps} \circ J^{-1}$ is càdlàg and $(\tau,v) \leq (\tau_n,v)$, by dominated convergence,
$\ov{u}^{\eps} \circ J^{-1}(\tau_n,v) \lra \ov{u}^{\eps} \circ J^{-1}(\tau,v)$ in $L^1(\Om, \calf, \bbp)$ as $\nto$ for any $0\leq v \leq 1$. Finally, by It\={o}'s isometry, the Lipschitz continuity of $f$ and with $\ind_{\emptyset} \equiv 0$,
\begin{equation} \label{NR13}
\begin{split}
& \bbe \left[ \left( u^{\eps} \circ J^{-1}(\tau_n,v) - u^{\eps} \circ J^{-1}(\tau,v) \right)^2 \right] \\ & \qquad \leq 
C \bbe \left[
\int_{\bbr^+ \times \bbr} \ind_{A^+(J^{-1}(\tau_n,v)) \setminus A^+(J^{-1}(\tau,v)) } (s,y) \left(\big \vert u^{\eps}(s,y) \big \vert^2 + 1 \right) \, \dd s \, \dd y
\right]
\end{split}
\end{equation}
%\\ & \qquad = 
%\bbe \left[ \left( \frac{1}{2}
%\int_{\bbr^+ \times \bbr} \ind_{A^+(J^{-1}(\tau_n,v)) \setminus A^+(J^{-1}(\tau,v)) } (s,y) \frac{f(u^{\eps}(s,y))}{\si(\eps)} \, L^{\eps}(\dd s, \dd y)
%\right)^2 \right] \\ 
for all $0 \leq v \leq 1$ and $n \in \bbn$. We infer $u^{\eps} \circ J^{-1}(\tau_n,v) \lra u^{\eps} \circ J^{-1}(\tau,v)$ in $L^2(\Om, \calf, \bbp)$ as $\nto$ for any $0\leq v \leq 1$, again by dominated convergence ($\ind_{A^+(J^{-1}(\tau_n,v)) \setminus A^+(J^{-1}(\tau,v)) } \lra 0$ pointwise on $\Om \times \bbr^+ \times \bbr$ as $\nto$ and the integrand above may be approximated by the integrable function $\ind_{A^+(u^*)}  (\vert u^{\eps} \vert^2 + 1)$ with $u^* = (3/2,1/2)$).
%We deduce $u^{\eps} \circ J^{-1}(\tau,v) = \ov{u}^{\eps} \circ J^{-1}(\tau,v)$ $\bbp$-almost surely for any $0\leq v \leq 1$ and $\eps > 0$.

Now we assume that each $u^{\eps}$ is separable and let $(\eps_n)_{n \in \bbn}$, $(h_n)_{n \in \bbn}$ be sequences of positive numbers with $\eps_n \lra 0$ and $h_n \lra 0$ as $n \rightarrow \infty$. Let also $(T_n)_{n \in \bbn}$ be a sequence of natural $1$-stopping times for $\ov{u}^{\eps_n} \circ J^{-1}$ with $T_n \in [0, 1]$. As for \eqref{NR13} and using \eqref{unif:bound:sol:Levy}, we obtain
\begin{equation} \label{NR14}
\begin{split}
&  \bbe \left[ \Big \vert
\ov{u}^{\eps_n} \circ J^{-1}(T_n + h_n, v) - \ov{u}^{\eps_n} \circ J^{-1}(T_n, v)
\Big \vert^2 \right] \! = \! \bbe \left[ \Big \vert
u^{\eps_n} \circ J^{-1}(T_n + h_n, v) - u^{\eps_n} \circ J^{-1}(T_n, v)
\Big \vert^2 \right] \\ & \qquad \leq 
C \bbe \left[ \left(\sup_{(s,y) \in A^+(u^*)} \big \vert u^{\eps_n}(s,y) \big \vert^2 + 1\right)
\int_{\bbr^+ \times \bbr} \ind_{A^+(J^{-1}(T_n + h_n,v)) \setminus A^+(J^{-1}(T_n,v)) } (s,y)  \, \dd s \, \dd y
\right] \\ & \qquad \leq 
C h_n \left(\bbe \left[ \sup_{(s,y) \in A^+(u^*)} \big \vert u^{\eps_n}(s,y) \big \vert^2 \right] + 1\right) \leq 
C h_n \left(\sup_{(s,y) \in A^+(u^*)} \bbe \left[ \vert u^{\eps_n}(s,y) \vert^2 \right] + 1\right) \leq C h_n,
\end{split}
\end{equation}
which goes to 0 as $\nto$ for all $0 \leq v \leq 1$. We used the inverse mapping of $H$ to see that whenever $T_n + v \geq 1$, the surface integral inside the third expectation in \eqref{NR14} equals
\begin{equation*}
\frac{9}{2} \left(T_n + h_n + v - \frac{1}{\sqrt{2}} \right)^2 - \frac{9}{2} \left(T_n + v - \frac{1}{\sqrt{2}} \right)^2 =
9 h_n \left( T_n + \frac{h_n}{2} + v - \frac{1}{\sqrt{2}} \right) (\leq C h_n).
\end{equation*}
Plus, the maximal inequality on the last line of \eqref{NR14} is a consequence of \eqref{max:in}. Analogously, if $(T_n)_{n \in \bbn}$ is a sequence of natural $2$-stopping times for $\ov{u}^{\eps_n} \circ J^{-1}$ with $T_n \in [0, 1]$,
\begin{equation*}
\sup_{0 \leq v \leq 1} \bbe \left[ \Big \vert \ov{u}^{\eps_n} \circ J^{-1}(v, T_n + h_n) - \ov{u}^{\eps_n} \circ J^{-1}(v, T_n) \Big \vert^2 \right]  \leq C h_n  \lra 0 \quad \textrm{as} \quad \nto.
\end{equation*}
Therefore, the random fields $\ov{u}^{\eps} \circ J^{-1}$ satisfy all conditions for Theorem~4.I in \cite{Ivanoff96} to apply. \qed
\epr

We now state a tightness result for $u^{\eps}$ and $\ov{u}^{\eps}$ in $L^2$-space.

\begin{Theorem} \label{theo:tight:L2:cadlag:vers:u:eps}
	The family $\{ u^{\eps} \mid \eps > 0 \}$ of mild solutions to \eqref{WAVE:Levy} is tight in the Hilbert space $L^2([0,T] \times I)$ for any $T > 0$ and finite interval $I \subseteq \bbr$.
\end{Theorem}

\bpr
Let $\{\Psi_{k} \mid k \in \bbn\}$ be a countable orthonormal basis of $L^2([0,T] \times I)$. By the stochastic Fubini theorem (see e.g. Theorem~2.6 in \cite{Walsh86}), for all $\eps > 0$ and $k \in \bbn$,
\begin{equation*}
\begin{split}
\langle u^{\eps}, \Psi_k \rangle & = \int_{0}^{T} \int_{I} \left( 
\int_{0}^{T} \int_{\bbr} G(t,x;s,y) \frac{f(u^{\eps}(s,y))}{\si(\eps)} \, L^{\eps}(\dd s, \textrm{d}y) 
\right) \Psi_k(t,x) \, \dd t \, \dd x \\ & = 
\int_{0}^{T} \int_{\bbr} \left( 
\int_{0}^{T} \int_{I} G(t,x;s,y) \Psi_k(t,x) \, \dd t \, \dd x
\right)  \frac{f(u^{\eps}(s,y))}{\si(\eps)} \, L^{\eps}(\dd s, \textrm{d}y) \quad \bbp\textrm{-almost surely.}
\end{split}
\end{equation*}
Using \eqref{unif:bound:sol:Levy}, Parseval's identity and Fubini's theorem, we infer
\begin{equation*}
\begin{split}
& \sum_{k=0}^{\infty} \sup_{\eps > 0} \bbe \left[ {\langle u^{\eps}, \Psi_k \rangle}^2 \right] \leq C 
\int_{0}^{T} \int_{\bbr} \left(\sum_{k=0}^{\infty} {\langle G( \cdot , \cdot ; s,y), \Psi_k \rangle}^2\right) \dd s \, \textrm{d}y \\ & \qquad = 
\int_{0}^{T} \int_{\bbr} \left(
\int_{0}^{T} \int_{I} G^2(t,x;s,y) \, \dd t \, \dd x
\right) \dd s \, \textrm{d}y  = 
\int_{0}^{T} \int_{I} \left(
\int_{0}^{T} \int_{\bbr} G^2(t,x;s,y) \, \dd s \, \dd y
\right) \dd t \, \dd x,
\end{split}
\end{equation*}
which is finite since the last inner integral equals $t^2/4$. This implies by Markov's inequality,
\begin{equation*}
\sup_{\eps > 0} \bbp \left( \sum_{k = N}^{\infty} {\langle u^{\eps}, \Psi_k \rangle}^2 > \de \right) \leq 
\frac{1}{\de} \sum_{k=N}^{\infty} \sup_{\eps > 0} \bbe \left[ {\langle u^{\eps}, \Psi_k \rangle}^2 \right] \lra 0 \quad \textrm{as} \quad N \to \infty
\end{equation*}
for all $\de > 0$ as well as 
\begin{equation*}
\sup_{\eps > 0} \bbp \left( \sum_{k=0}^{N} {\langle u^{\eps}, \Psi_k \rangle}^2 > \de \right) \leq 
\frac{1}{\de} \sum_{k=0}^{\infty} \sup_{\eps > 0} \bbe \left[ {\langle u^{\eps}, \Psi_k \rangle}^2 \right] \lra 0 \quad \textrm{as} \quad \de \to \infty
\end{equation*}
for all $N \in \bbn$. So we can apply Theorem 1 in \cite{Suquet} and conclude the proof. \qed
\epr

We turn to the $H_{-r}(\bbr)$-valued process $v^{\eps}$ in \eqref{dist:proc:v} and first show that it has a càdlàg version.

\bpr[of Theorem \ref{theo:cadlag:vers:v:eps}]
The proof relies on the Hilbert space structure of $H_{-r}(\bbr)$. First, we show that for any $\phi \in \cals(\bbr)$ and $\eps > 0$, the real-valued process $(\langle v^{\eps}_t, \phi \rangle)_{t \geq 0}$ has a càdlàg modification. Use \eqref{dist:deriv:G}, \eqref{dist:proc:v} and the fundamental theorem of calculus to rewrite for all $t \geq 0$, 
\begin{equation*}
	\begin{split}
		\langle v^{\eps}_t, \phi \rangle & =
		\int_{0}^{t} \int_{\bbr} \phi(y) \frac{f(u^{\eps}(s,y))}{\si(\eps)} \, L^{\eps}(\dd s, \dd y) \\ & \quad \,\, +
		\frac{1}{2} \int_{0}^{t} \int_{\bbr} \left( \int_{s}^{t} \phi'(y+(r-s)) \, \dd r - \int_{s}^{t} \phi'(y-(r-s)) \, \dd r \right)
		\frac{f(u^{\eps}(s,y))}{\si(\eps)} \, L^{\eps}(\dd s, \dd y)
	\end{split}
\end{equation*}
%= \frac{1}{2}
%\int_{0}^{t} \int_{\bbr} \phi(y \pm (t-s)) \frac{f(u^{\eps}(s,y))}{\si(\eps)} \, L^{\eps}(\dd s, \dd y) \\ & =
$\bbp$-almost surely, and then the stochastic Fubini theorem on the last double integral to obtain
\begin{equation} \label{semimart:repr}
	\langle v^{\eps}_t, \phi \rangle = 
	\int_{0}^{t} \int_{\bbr} \phi(y) \frac{f(u^{\eps}(s,y))}{\si(\eps)} \, L^{\eps}(\dd s, \dd y) +
	\frac{1}{2} \int_{0}^{t} J_r^{\eps}(\phi) \, \dd r 
	\quad \bbp \textrm{-almost surely}
\end{equation}
where we have set
\begin{equation*}
	J_r^{\eps}(\phi) = \int_{0}^{r} \int_{\bbr}  \left(\phi'(y+(r-s)) -  \phi'(y-(r-s))\right) \frac{f(u^{\eps}(s,y))}{\si(\eps)} \, L^{\eps}(\dd s, \dd y), \quad r \geq 0.
\end{equation*}
The semimartingale on the right-hand side of \eqref{semimart:repr}, that we will denote by $X^{\eps}(\phi)$, is càdlàg. 

Next, fix an arbitrary $T > 0$ and $\eps > 0$. Doob's inequality, It\={o}'s isometry and \eqref{Bound:sol:Levy} yield
\begin{equation} \label{NR5}
	\begin{split}
		\bbe \left[ \sup_{t \leq T} \Big \vert X_t^{\eps}(\phi) \Big \vert^2 \right] \leq
		C \int_{0}^{T} \int_{\bbr} \phi^2(y)  \, \dd s \, \dd y + C \int_{0}^{T} 
		\left( \int_{0}^{r} \int_{\bbr}  \phi'(y \pm (r-s))^2 \, \dd s \, \dd y \right) \dd r
	\end{split}
\end{equation}
%& \leq 
%C \bbe \left[ \left(\int_{0}^{T} \int_{\bbr} \phi(y) \frac{f(u^{\eps}(s,y))}{\si(\eps)} \, L^{\eps}(\dd s, \dd y)\right)^2 \right] + 
%C  \int_{0}^{T} \bbe \left[ \vert J_r^{\eps}(\phi) \vert^2 \right] \dd r \\ &  
for all $\phi \in \cals(\bbr)$. Now the Hermite functions $h_q$ satisfy the recursion relation $h_q'(x) = \sqrt{q/2} h_{q-1}(x) - \sqrt{(q+1)/2} h_{q+1}(x)$ for all $q \in \bbn$ and $x \in \bbr$, from which we obtain by orthogonality,
\begin{equation} \label{norm:deriv:Herm}
	\int_{\bbr} h_q'(x)^2 \, \dd x = \frac{q}{2} \int_{\bbr} h^2_{q-1}(x) \dd x + \frac{q + 1}{2} \int_{\bbr} h^2_{q+1}(x) \dd x = q + \frac{1}{2}.
\end{equation}
We carry forward the estimation in \eqref{NR5} for $\phi = h_q$, whence
\begin{equation} \label{unif:bound:semi}
	\bbe \left[ \sup_{t \leq T} \Big \vert X_t^{\eps}(h_q) \Big \vert^2 \right] \leq 
	C \left(1 + \int_{\bbr} h_q'(x)^2 \, \dd x \right) \leq C(1 + 2q) \quad \textrm{for all} \quad q \in \bbn.
\end{equation}
Fix $r > 2$. It is easy to see that for each $N \in \bbn$, the $H_{-r}(\bbr)$-valued process 
\begin{equation} \label{cadlag:proc:N}
	\sum_{q=0}^{N} (1 + 2q)^{-r/2} X^{\eps}(h_q) e_{q,-r}
\end{equation}
with $e_{q,-r}$ as in \eqref{ONB:H-r}, is càdlàg. Recall the Fourier expansion \eqref{dual:ONB:repr} and use \eqref{unif:bound:semi} to obtain
\begin{equation} \label{NR8}
	\begin{split}
		\bbe \left[\sup_{t \leq T} {\left \Vert \sum_{q=N+1}^{M} (1 + 2q)^{-r/2} X_t^{\eps}(h_q) e_{q,-r} \right \Vert}_{-r}^2\right] &  \leq 
		\sum_{q=N+1}^{M} (1 + 2q)^{-r} \bbe \left[\sup_{t \leq T} \Big \vert X_t^{\eps}(h_q) \Big \vert^2 \right] \\ & \leq
		C \sum_{q=N+1}^{M} (1 + 2q)^{-r+1} \lra 0
	\end{split}
\end{equation}
as $N, M \to \infty$ since $r > 2$. 
%Without loss of generality, this convergence has the following interpretation: for each $\eps > 0$, the processes in \eqref{cadlag:proc:N} build $\bbp$-almost surely a Cauchy sequence in the complete metric space $(D([0,T], H_{-r}(\bbr)), \sup_{t \leq T} {\Vert \cdot \Vert}_{-r})$. 
Consequently, standard arguments show that there exists a process $\ov{v}^{\eps} \in D([0,T], H_{-r}(\bbr))$ such that $\bbp$-almost surely,
\begin{equation} \label{cadlag:vers:v:eps} 
\ov{v}^{\eps}_t = \sum_{q=0}^{\infty} (1 + 2q)^{-r/2} X_t^{\eps}(h_q) e_{q,-r} \quad \textrm{in} \quad H_{-r}(\bbr) \quad \textrm{for all} \quad t \leq T.
\end{equation}
By \eqref{semimart:repr}, this process is a version of $(v^{\eps}_t)_{t \leq T}$ in $H_{-r}(\bbr)$.
% since for all $t \leq T$,
%\begin{equation*}
%	\bbe \left[ {\Vert \ov{v}_t^{\eps} - v^{\eps}_t \Vert}_{-r}^2 \right] = 
%	\sum_{q=0}^{\infty} (1 + 2q)^{-r} \bbe \left[ \left( X_t^{\eps}(h_q) - \langle v^{\eps}_t, h_q \rangle \right)^2 \right] = 0.
%\end{equation*}
\qed
% The existence of a càdlàg version of $v^{\eps}$ on $\bbr^+$ for each $\eps > 0$ immediately follows.  \qed
\epr

Next, we show tightness of the càdlàg version $\ov{v}^{\eps}$.

\begin{Theorem} \label{theo:tight:Skor:cadlag:vers:v:eps}
	The family of processes $\{ \ov{v}^{\eps} \mid \eps > 0 \}$ where $\ov{v}^{\eps}$ is the càdlàg version of $v^{\eps}$ obtained in Theorem~\ref{theo:cadlag:vers:v:eps}, is tight in the Skorokhod space $D([0,T], H_{-r}(\bbr))$ for any $r > 2$ and $T > 0$. 
\end{Theorem}

\bpr
We first check that $\{ \ov{v}^{\eps} \mid \eps > 0 \}$ satisfies the Aldous condition for tightness. To this end, let $(\varepsilon_n)_{n \in \mathbb{N}}$ and $(h_n)_{n \in \mathbb{N}}$ be sequences of positive numbers with $\varepsilon_n \rightarrow 0$ and $h_n \rightarrow 0$ as $n \rightarrow \infty$. In addition, for each $n \in \mathbb{N}$, let $\tau_n  \in [0,T]$ be a stopping time with respect to the filtration generated by the process $(\ov{v}_t^{\eps_n})_{t \leq T}$. We will show
\begin{equation} \label{Aldous:cond}
	\mathbb{E}\left[
	{\Vert \ov{v}^{\eps_n}_{\tau_n + h_n} - \ov{v}^{\eps_n}_{\tau_n} \Vert}_{-r}^2
	\right] \lra 0 \quad \textrm{as} \quad  n \to \infty.
\end{equation}
Recall the series representation \eqref{cadlag:vers:v:eps} of $\ov{v}^{\eps_n}$, where $X^{\eps_n}(h_q)$ is the right-hand side of \eqref{semimart:repr} with $\phi = h_q$. We have for each $q \in \bbn$,
\begin{equation} \label{SM:repr:Aldous}
	\begin{split}
		X_{\tau_n + h_n}^{\eps_n}(h_q) - X_{\tau_n}^{\eps_n}(h_q)  & = 
		\int_{\tau_n}^{\tau_n + h_n} \int_{\bbr} h_q(y) \frac{f(u^{\eps_n}(s,y))}{\si(\eps_n)} \, L^{\eps_n}(\dd s, \dd y) +
		\frac{1}{2} \int_{\tau_n}^{\tau_n + h_n} J_r^{\eps_n}(h_q) \, \dd r \\ & 
		=: I_{q,n} + J_{q,n}.
	\end{split}
\end{equation}
We estimate the second moment of each of these two terms. For the first one, by It\={o}'s isometry and the Lipschitz continuity of $f$,
\begin{equation} \label{NR6}
	\begin{split}
		\bbe \left[ I_{q,n}^2 \right] & = 
		\bbe \left[ \int_{0}^{T} \int_{\bbr}  h^2_q(y) \ind_{(\tau_n, \tau_n + h_n]}(s)  f^2(u^{\eps_n}(s,y)) \, \dd s \, \dd y  \right] \\ & \leq 
		C \bbe \left[ \int_{0}^{T} \int_{\bbr}  h^2_q(y) \ind_{(\tau_n, \tau_n + h_n]}(s) \vert u^{\eps_n}(s,y) \vert^2 \, \dd s \, \dd y \right] +
		C \bbe \left[ \int_{0}^{T} \ind_{(\tau_n, \tau_n + h_n]}(s) \, \dd s \right] \\ & =
		C \int_{\bbr}  h^2_q(y) \bbe \left[ \int_{0}^{T} \ind_{(\tau_n, \tau_n + h_n]}(s) \vert u^{\eps_n}(s,y) \vert^2 \, \dd s \right] \dd y + 
		C h_n.
	\end{split}
\end{equation}
Furthermore, by the maximal inequality \eqref{max:in} (assuming separability),
\begin{equation} \label{Cairoli}
	\begin{split}
		\bbe \left[ \sup_{(s,y) \in [(0,x), (T,x)]_{\lpo} } \big \vert u^{\eps}(s,y) \big \vert^2  \right] \leq 
		\bbe \left[ \big \vert u^{\eps}(T,x) \big \vert^2 \right] \quad \textrm{for all} \quad x \in \bbr \quad \textrm{and} \quad \eps > 0.
	\end{split}
\end{equation}
Hence, the remaining integral on the right-hand side of \eqref{NR6} can further be estimated by
\begin{equation} \label{NR19}
	\begin{split}
		& \int_{\bbr}  h^2_q(y) \bbe \left[ \sup_{(s,z) \in [(0,y), (T,y)]_{\lpo} } \big \vert u^{\eps_n}(s,z) \big \vert^2 \, 
		\int_{0}^{T} \ind_{(\tau_n, \tau_n + h_n]}(s) \, \dd s \right] \dd y \\ & \qquad \quad =
		h_n \int_{\bbr}  h^2_q(y) \bbe \left[ \sup_{(s,z) \in [(0,y), (T,y)]_{\lpo} } \big \vert u^{\eps_n}(s,z) \big \vert^2 \right] \dd y \leq 
		h_n \int_{\bbr}  h^2_q(y) \bbe \left[ \big \vert u^{\eps_n}(T,y) \big \vert^2 \right] \dd y \\ & \qquad \quad 
		\leq C h_n \int_{\bbr}  h^2_q(y) \dd y = C h_n \quad \textrm{for all} \quad  q, n \in \bbn,
	\end{split}
\end{equation}
where \eqref{unif:bound:sol:Levy} was used for the last inequality. Note a significant difference here with the stochastic heat equation addressed in \cite{CC:TD}: The mild solution to that equation is not a multiparameter martingale, so instead of maximal inequalities as \eqref{Cairoli}, the factorization method from \cite{DaPrato87, Sanz} was used to prove the Aldous condition, see in particular Lemma~3.3 and (3.13) in \cite{CC:TD}.  

Next, by the same calculations as in \eqref{NR5} (but using \eqref{unif:bound:sol:Levy} instead of \eqref{Bound:sol:Levy}) and \eqref{unif:bound:semi}, 
$\sup_{\eps > 0} \bbe \left[ \vert J_r^{\eps}(h_q) \vert^2 \right] \leq C(1 + 2q)$ for all $q \in \bbn$ and $r \leq T$, so using Hölder's inequality,
\begin{equation*}
	\bbe \left[ J_{q,n}^2 \right] = \frac{1}{4} \bbe \left[ \left( \int_{0}^{T} \ind_{(\tau_n, \tau_n + h_n]} J_r^{\eps_n}(h_q) \, \dd r \right)^2 \right] 
	\leq C \bbe \left[ h_n \int_{0}^{T} \vert J_r^{\eps_n}(h_q) \vert^2 \, \dd r \right] 	
	\leq C h_n (1 + 2q)
\end{equation*}
for all $n, q \in \bbn$. Combine this with \eqref{SM:repr:Aldous}, \eqref{NR6} and \eqref{NR19} to obtain altogether
\begin{equation*}
	\begin{split}
		\mathbb{E}\left[
		{\Vert \ov{v}^{\eps_n}_{\tau_n + h_n} - \ov{v}^{\eps_n}_{\tau_n} \Vert}_{-r}^2
		\right] & = 
		\sum_{q=0}^{\infty} (1 + 2q)^{-r} \mathbb{E}\left[ \left(X_{\tau_n + h_n}^{\eps_n}(h_q) - X_{\tau_n}^{\eps_n}(h_q)\right)^2 \right] \\ & \leq
		2 \sum_{q=0}^{\infty} (1 + 2q)^{-r} \left(\mathbb{E}\left[ I_{q,n}^2 \right] + \mathbb{E}\left[ J_{q,n}^2 \right]\right) \leq
		h_n C \sum_{q=0}^{\infty} (1 + 2q)^{-r + 1} \lra 0
	\end{split}
\end{equation*}
as $\nto$ since $r > 2$, which is \eqref{Aldous:cond}.

In addition,
$\sum_{q=0}^{\infty} (1 + 2q)^{-r} \sup_{\eps > 0} \bbe [ X_{t}^{\eps}(h_q)^2 ] < \infty$  by \eqref{unif:bound:sol:Levy}, \eqref{unif:bound:semi} and since $r > 2$, so we can readily deduce, as in the proof of Theorem~\ref{theo:tight:L2:cadlag:vers:u:eps}, that the random elements $\{ \ov{v}_t^{\eps} \mid \eps > 0 \}$ are tight in $H_{-r}(\bbr)$ for any fixed $t \leq T$.

The claim of the theorem now directly follows from Theorem~6.8 in \cite{Walsh86}. \qed
%the random elements $\{ \ov{v}_t^{\eps} \mid \eps > 0 \}$ are tight in $H_{-r}(\bbr)$ for all fixed $t \leq T$. We have indeed by \eqref{unif:bound:semi} and since $r > 2$,
%\begin{equation*}
%	\sum_{q=0}^{\infty} (1 + 2q)^{-r} \sup_{\eps > 0} \bbe \left[ X_{t}^{\eps}(h_q)^2 \right] \leq 
%	\sum_{q=0}^{\infty} (1 + 2q)^{-r} \sup_{\eps > 0} \bbe \left[ \sup_{t \leq T} \Big \vert X_t^{\eps}(h_q) \Big \vert^2 \right] \leq 
%	C \sum_{q=0}^{\infty} (1 + 2q)^{-r+1} < \infty.
%\end{equation*}
%This implies by Markov's inequality and \eqref{cadlag:vers:v:eps},
%\begin{equation*}
%	\sup_{\eps > 0} \bbp \left( \sum_{q = N}^{\infty} \langle \ov{v}_t^{\eps}, e_{q,-r} \rangle^2_{-r} > \de \right) \leq 
%	\frac{1}{\de} \sum_{q=N}^{\infty} (1 + 2q)^{-r} \sup_{\eps > 0} \bbe \left[ X_{t}^{\eps}(h_q)^2 \right] \lra 0 \,\,\, \textrm{as} \,\,\, N \to \infty
%\end{equation*}
%for all $\de > 0$ as well as 
%\begin{equation*}
%	\sup_{\eps > 0} \bbp \left( \sum_{q =0}^{N} \langle \ov{v}_t^{\eps}, e_{q,-r} \rangle^2_{-r} > \de \right) \leq 
%	\frac{1}{\de} \sum_{q=0}^{\infty} (1 + 2q)^{-r} \sup_{\eps > 0} \bbe \left[ X_{t}^{\eps}(h_q)^2 \right] \lra 0 \,\,\, \textrm{as} \,\,\, \de \to \infty
%\end{equation*}
%for all $N \in \bbn$. So we can apply Theorem 1 in \cite{Suquet}.
\epr

We end this section with a tightness result for $v^{\eps}$ and $\ov{v}^{\eps}$ in $L^2$-space.

\begin{Theorem} \label{tight:v:eps:inL2}
	The distribution-valued processes $\{ v^{\eps} \mid \eps > 0 \}$ with $v^{\eps}$ as in \eqref{dist:proc:v}, are tight in the Hilbert space $L^2([0,T], H_{-r}(\bbr))$ for each $r > 1$ and $T > 0$.
\end{Theorem}

\bpr
First, each $v^{\eps}$ is an element of $L^2([0,T], H_{-r}(\bbr))$ as is seen from
\begin{equation*}
	\bbe \left[ \int_{0}^{T} {\Vert v^{\eps}_t \Vert}_{-r}^2 \, \dd t \right] = 
	\sum_{q=0}^{\infty} (1 + 2q)^{-r} \int_{0}^{T} \bbe \left[ {\langle v_t^{\eps}, h_q \rangle}^2 \right] \dd t \leq C \sum_{q=0}^{\infty} (1 + 2q)^{-r} < \infty
\end{equation*}
which follows from \eqref{NR4} and $r > 1$. 

The scalar product in $L^2([0,T], H_{-r}(\bbr))$ is given by
$\langle f, g \rangle = \int_{0}^{T} {\langle f_t, g_t \rangle}_{-r} \, \dd t $ and it is easy to see that an orthonormal basis is formed by $\{ \phi_i e_{q,-r} \mid i, q \in \bbn \}$ with $\phi_i(t) = \sqrt{2/T} \sin(it\pi /T)$ and $e_{q,-r}$ as in \eqref{ONB:H-r}.
%Indeed, we have
%\begin{equation*}
%	\begin{split}
%		\langle \phi_i e_{q,-r}, \phi_j e_{p,-r} \rangle & = \left(\int_{0}^{T} \phi_i(t) \phi_j(t)  \, \dd t\right) {\langle e_{q,-r}, e_{p,-r} \rangle}_{-r} = \de_{ij} \de_{qp} \quad \textrm{as well as} \\ 
%		{\Vert f \Vert}^2 & = \int_{0}^{T} \left( \sum_{q=0}^{\infty} {\langle f_t, e_{q,-r} \rangle}_{-r}^2 \right) \dd t = 
%		\sum_{q=0}^{\infty} \int_{0}^{T}  {\langle f_t, e_{q,-r} \rangle}_{-r}^2 \dd t \\ & = 
%		\sum_{q=0}^{\infty} \sum_{i=1}^{\infty} \left(\int_{0}^{T} {\langle f_t, e_{q,-r} \rangle}_{-r} \phi_i(t) \, \dd t \right)^2 
%		= \sum_{i, q = 0}^{\infty} {\langle f, \phi_i e_{q,-r} \rangle}^2.
%	\end{split}
%\end{equation*}
By \eqref{dist:proc:v} and the stochastic Fubini theorem, for all $i, q \in \bbn$ and $\eps > 0$,
\begin{equation*}
	\begin{split}
		\int_{0}^{T} {\langle v^{\eps}_t, h_q \rangle} \phi_i(t) \, \dd t  = \frac{1}{2}
		\int_{0}^{T} \int_{\bbr} \left( \int_{s}^{T}  h_q(y \pm (t-s)) \phi_i(t) \, \dd t \right) \frac{f(u^{\eps}(s,y))}{\si(\eps)} \, L^{\eps}(\dd s, \dd y) 
		\quad \bbp \textrm{-a.s.}
	\end{split}
\end{equation*}
%& = \frac{1}{2}
%\int_{0}^{T} \left( \int_{0}^{t} \int_{\bbr} h_q(y \pm (t-s)) \frac{f(u^{\eps}(s,y))}{\si(\eps)} \, L^{\eps}(\dd s, \dd y) \right) \phi_i(t) \, \dd t \\ &
Therefore, by duality, It\={o}'s isometry and \eqref{unif:bound:sol:Levy}, we have
\begin{equation*}
	\begin{split}
		\bbe \left[ {\langle v^{\eps}, \phi_i e_{q,-r} \rangle}^2 \right] \leq 
		C (1 + 2q)^{-r} \int_{0}^{T} \int_{\bbr} \left( \int_{s}^{T}  h_q(y \pm (t-s)) \phi_i(t) \, \dd t \right)^2  \dd s \, \dd y.
	\end{split}
\end{equation*} 
% & = (1 + 2q)^{-r} \bbe \left[ \left(\int_{0}^{T} {\langle v^{\eps}_t, h_q \rangle} \phi_i(t) \, \dd t\right)^2 \right] \\ &
Using Parseval's identity relative to the orthonormal basis of $L^2([0,T])$, we obtain altogether
\begin{equation} \label{NR9}
	\begin{split}
		\sum_{i, q = 0}^{\infty} \sup_{\eps > 0} \bbe \left[ {\langle v^{\eps}, \phi_i e_{q,-r} \rangle}^2 \right] & \leq 
		C \sum_{q=0}^{\infty} (1 + 2q)^{-r} \int_{0}^{T} \int_{\bbr} \, \sum_{i=0}^{\infty} \left( \int_{s}^{T}  h_q(y \pm (t-s)) \phi_i(t) \, \dd t \right)^2  \dd s \, \dd y \\ & =
		C \! \sum_{q=0}^{\infty} (1 + 2q)^{-r} \! \int_{0}^{T} \int_{\bbr} \int_{s}^{T}  h^2_q(y \pm (t-s))  \, \dd t \, \dd y \, \dd s \leq C \sum_{q=0}^{\infty} (1 + 2q)^{-r}
	\end{split}
\end{equation}
which is finite since $r > 1$. We can now conclude analogously to the proof of Theorem~\ref{theo:cadlag:vers:v:eps}. \qed
\epr

\begin{Corollary} \label{cor:tight:L2:v:eps}
	The distribution-valued processes $\{ \ov{v}^{\eps} \mid \eps > 0 \}$ where $\ov{v}^{\eps}$ is the càdlàg version of $v^{\eps}$ in Theorem~\ref{theo:cadlag:vers:v:eps}, are tight in $L^2([0,T], H_{-r}(\bbr))$ for any $r > 2$ and $T > 0$.
\end{Corollary}

\bpr Since $\ov{v}_t^{\eps}$ lives in $H_{-r}(\bbr)$ for $r > 2$ only, this is a direct consequence of Theorem~\ref{tight:v:eps:inL2}. \qed
\epr

\subsection{Proofs for Section \ref{sec:main:res}}

We first show that a mild solution to \eqref{WAVE:Levy} satisfies equation \eqref{weak:form:formal} with $\pd_t u^{\eps}$ replaced by $v^{\eps}$.

\begin{Proposition} \label{prop:mild:implies:weak}
	Let $u^{\eps}$ be a mild solution to \eqref{WAVE:Levy} and $v^{\eps}$ the process defined in \eqref{dist:proc:v}. For each $\eps > 0$, the pair $(u^{\eps}, v^{\eps})$ satisfies the following weak formulation of the stochastic wave equation on $\bbr^+ \times \bbr$: For any $\phi_1, \phi_2 \in C_c^{\infty}(\bbr)$ and $t \geq 0$, we have 
	\begin{equation} \label{weak:form:exact:Levy}
	\begin{split}
	& \int_{\bbr} u^{\eps}(t,x) \phi_1(x) \, \dd x + \langle v_t^{\eps}, \phi_2 \rangle \\ & \qquad \quad = 
	\int_{0}^{t} \left(\int_{\bbr} u^{\eps}(s,x) \phi_2''(x) \, \dd x +  \langle v_s^{\eps}, \phi_1 \rangle \right) \dd s + 
	\int_{0}^{t} \int_{\bbr} \phi_2(x) \disfrac{f(u^{\eps}(s,x))}{\si(\eps)} \, L^{\eps}(\dd s, \dd x) \quad \bbp \textrm{-a.s}.
	\end{split}
	\end{equation}
\end{Proposition}

\bpr 
Recall first \eqref{dist:deriv:G} and apply the stochastic Fubini theorem in order to obtain
\begin{equation} \label{NR21}
\begin{split}
\int_{0}^{t} \int_{\bbr} u^{\eps}(s,x) \phi_2''(x) \, \dd x \, \dd s & = 
\int_{0}^{t} \int_{\bbr} \left( \int_{0}^{t} \int_{\bbr} G(s,x; r,y) \phi_2''(x) \, \dd s \, \dd x \right) \disfrac{f(u^{\eps}(r,y))}{\si(\eps)} \, L^{\eps}(\dd r, \dd y) 
\quad \textrm{and} \\
\int_{0}^{t} \langle v_s^{\eps}, \phi_1 \rangle \, \dd s & = 
\int_{0}^{t} \int_{\bbr} \left( \int_{0}^{t} \int_{\bbr} \phi_1(x) \frac{\dd  G}{\dd x}(s, \dd x; r, y) \, \dd s \right)  \disfrac{f(u^{\eps}(r,y))}{\si(\eps)} \, L^{\eps}(\dd r, \dd y).
\end{split}
\end{equation}
We further calculate for both inner integrals in \eqref{NR21} and fixed $0 \leq r \leq t$ and $y \in \bbr$,
\begin{equation*}
\begin{split}
\int_{0}^{t} \int_{\bbr}  G(s,x;r,y) \phi_2''(x) \, \dd x \, \dd s & = 
\frac{1}{2} \phi_2(y \pm (t-r)) - \phi_2(y) = \int_{\bbr} \phi_2(x)  
\frac{\dd  G}{\dd x}(t, \dd x; r, y) - \phi_2(y) \quad \textrm{and} \\
\int_{0}^{t} \int_{\bbr} \phi_1(x) \frac{\dd  G}{\dd x}(s, \dd x; r, y) \, \dd s & = \int_{\bbr} G(t,z;r,y) \phi_1(z) \, \dd z.
\end{split}
\end{equation*}
Now insert the last integral accordingly in \eqref{NR21} and apply again the stochastic Fubini theorem. \qed
\epr

The next theorem is a converse of Proposition~\ref{prop:mild:implies:weak} in the following sense: If a random field on $[0,T] \times [-T,L+T]$ satisfies (together with an auxiliary distribution-valued process) the weak formulation of the stochastic wave equation (on $\bbr^+ \times \bbr$) driven by Gaussian noise "restricted" to $[0,T] \times [-T,L+T]$, 
%(in particular, the test functions now considered must have support in $(-T,L+T)$)
then it is a mild solution to \eqref{WAVE:Gauss} on $[0,T] \times [0,L]$.

\begin{Theorem} \label{theo:weak:implies:mild}	
	On a complete stochastic basis $(\wt{\Om}, \wt{\calf}, \wt{\bfil}, \wt{\bbp})$, let $\wt{W}$ be a Gaussian space--time white noise on $[0,T] \times [-T,L+T]$ for some $T > 0$ and $L > 0$. Assume we have a $\lpo$-càdlàg random field $w=\{ w(t,x) \mid (t,x) \in [0,T] \times [-T,L+T] \}$ satisfying
	\begin{equation} \label{ess:sup:2}
	\esssup_{(t,x) \in [0,T] \times [-T,L+T]} \bbe \left[ \big \vert w(t,x) \big \vert^2 \right] < \infty
	\end{equation}
	and a $H_{-r}(\bbr)$-valued càdlàg process $(\theta_t)_{t \leq T}$ for some $r > 2$. Assume for any $x \in [-T,L+T]$, $w(0,x)=\theta_0=0$ $\wt{\bbp}$-a.s. and that for all $\phi_1, \phi_2 \in C_c^{\infty}((-T,L+T))$, the pair $(w, \theta)$ satisfies \eqref{weak:form:exact:Gauss} with probability one. Then $w$ is on $[0,T] \times [0,L]$ the continuous mild solution to the stochastic wave equation \eqref{WAVE:Gauss:2} driven by $\dot{\wt{W}}$, and $\theta$ satisfies \eqref{NR18} for all $t \leq T$ and $\phi \in C_c^{\infty}((0,L))$.
\end{Theorem}

For the proof of this theorem, we need the following technical lemma.

\begin{Lemma} \label{lem:dense}
	Let $T > 0$ and $I \subseteq \bbr$ be a finite \emph{open} interval. The tensor product $C^{\infty}([0,T]) \otimes C_c^{\infty}(I)$ is dense in $C_c^{\infty}([0,T] \times I)$ with respect to each norm $\sum_{\vert \al \vert \leq N} {\Vert \cdot \Vert}_{\infty, \al}$ with $N \in \bbn$ and
	\begin{equation} \label{seminorm}
	{\Vert f \Vert}_{\infty, \al} = {\Vert f^{(\al)} \Vert}_{\infty} = \sup_{(t,x) \in [0,T] \times I} \Big \vert f^{(\al)}(t,x) \Big \vert
	\end{equation}
	and $f^{(\al)} = \pd_t^{\al_1} \pd_x^{\al_2} f$ with multi-index $\al = (\al_1, \al_2) \in \bbn^2$.
\end{Lemma}

\bpr
Assume for simplicity $I = (0, 1)$, fix $f \in C_c^{\infty}([0,T] \times I)$ and a compact set $A \subseteq I$ such that $\supp f \subseteq [0,T] \times A$. Furthermore, let $b$ be a $C_c^{\infty}(\bbr)$-function such that $0 \leq b \leq 1$, $b \equiv 1$ on $A$ and $A \subseteq \supp b \subsetneq I$. Set $K = \supp b$.

The set of all polynomials on $[0,T] \times K$ is dense in $C^{\infty}([0,T] \times K)$ with respect to each norm $\sum_{\vert \al \vert \leq N} {\Vert \cdot \Vert}_{\infty, \al}$ (with the obvious restriction of domain of definition). We prove this by induction on the differentiation order $N$: If $N=0$, it is a direct consequence of the Stone--Weierstrass theorem and if the claim holds for $N-1$, choose $g \in C^{\infty}([0,T] \times K)$ and write $g(t,x) = \int_{0}^{t} \pd_t g(s,x) \, \dd s + \int_{a}^{x} \pd_x g(0,y) \, \dd y + g(a,0)$ (assuming $K=[a,b]$ for simplicity). By assumption, we can find polynomials $A_n$, resp. $B_n$, that converge in $\sum_{\vert \al \vert \leq N - 1} {\Vert \cdot \Vert}_{\infty, \al}$ to $\pd_{t} g$, resp. $\pd_{x} g$. Then the polynomial $C_n(t,x) = \int_{0}^{t} \pd_t A_n(s,x) \, \dd s + \int_{a}^{x} B_n(0,y) \, \dd y + g(a,0)$ converges to $g$ in $\sum_{\vert \al \vert \leq N} {\Vert \cdot \Vert}_{\infty, \al}$.
%the case $C([0,T] \times K)$ is a direct consequence of the Stone--Weierstrass theorem, and if we assume that the claim holds for $C^{k-1}([0,T] \times K)$, take any $g \in C^{k}([0,T] \times K)$ and write $g(t,x) = \int_{0}^{t} \pd_t g(s,x) \, \dd s + \int_{0}^{x} \pd_x g(0,y) \, \dd y + g(0,0)$. By assumption, find polynomials $A_n$, resp. $B_n$, that converge to $\pd_{t} g$, resp. $\pd_{x} g$, in each ${\Vert \cdot \Vert}_{\infty, \al}$ with $\vert \al \vert \leq k-1$. Then the polynomial $C_n(t,x) = \int_{0}^{t} \pd_t A_n(s,x) \, \dd s + \int_{0}^{x} B_n(0,y) \, \dd y + g(0,0)$ converges to $g$ in each ${\Vert \cdot \Vert}_{\infty, \al}$ with $\vert \al \vert \leq k$.

Now fix $N \in \bbn$ and choose a sequence of polynomials
$P_n(t,x) = \sum_{i,j = 1}^{N_n} \al_{i,n} \beta_{j,n} t^i x^j$ with $\al_{i,n}, \beta_{j,n} \in \bbr$ and $N_n \in \bbn$ such that $P_n^{(\al)}$ converges to $f^{(\al)}$ uniformly on $[0,T] \times K$ for all $\vert \al \vert \leq N$. We can write $b(x) P_n(t,x) = \sum_{i, j=0}^{N_n} \al_{i,n} \beta_{j,n} t^i (b(x) x^j)$, so each $b P_n$ lies in $C^{\infty}([0,T]) \otimes C_c^{\infty}(I)$ since $b \in C_c^{\infty}(I)$. We now make the following calculations on $[0,T] \times I$. By the Leibniz rule,
${\Vert (f - bP_n)^{(0,k)} \Vert}_{\infty} \leq 
{\Vert f^{(0,k)} - b P_n^{(0,k)} \Vert}_{\infty} + C \sum_{l=1}^{k} {\Vert b^{(l)} P_n^{(l,k-l)} \Vert}_{\infty}$ for any $k \in \bbn$. The first term on the right-hand side can further be estimated by
\begin{equation*}
\sup_{(t,x) \in [0,T] \times A} \Big \vert f^{(0,k)}(t,x) - P_n^{(0,k)}(t,x) \Big \vert + 
\sup_{(t,x) \in [0,T] \times ( K \setminus A)} \Big \vert P_n^{(0,k)}(t,x) \Big \vert,
\end{equation*}
%\begin{equation*}
%\begin{split}
%& {\Vert f^{(0,k)} - b P_n^{(0,k)} \Vert}_{\infty} \\ & \qquad \leq 
%\sup_{(t,x) \in [0,T] \times A} \Big \vert f^{(0,k)}(t,x) - b(x) P_n^{(0,k)}(t,x) \Big \vert + 
%\sup_{(t,x) \in [0,T] \times ( I \setminus A)} \Big \vert f^{(0,k)}(t,x) - b(x) P_n^{(0,k)}(t,x) \Big \vert \\ & \qquad \leq
%\sup_{(t,x) \in [0,T] \times A} \Big \vert f^{(0,k)}(t,x) - P_n^{(0,k)}(t,x) \Big \vert + 
%\sup_{(t,x) \in [0,T] \times ( K \setminus A)} \Big \vert P_n^{(0,k)}(t,x) \Big \vert,
%\end{split}
%\end{equation*}
which goes to 0 for each $k \leq N$ by assumption on $P_n$. On the other hand,
\begin{equation*}
\sum_{l=1}^{k} {\Vert b^{(l)} P_n^{(l,k-l)} \Vert}_{\infty}  \leq 
\left(\max_{l=1, \dots, k} {\Vert b^{(l)} \Vert}_{\infty}\right) \sum_{l=1}^{k} \sup_{(t,x) \in [0,T] \times (K \setminus A)} \Big \vert P_n^{(l,k-l)}(t,x) \Big \vert,
\end{equation*}
%\begin{equation*}
%\begin{split}
%\sum_{l=1}^{k} {\Vert b^{(l)} P_n^{(l,k-l)} \Vert}_{\infty} & = 
%\sum_{l=1}^{k} \sup_{(t,x) \in [0,T] \times (K \setminus A)} \Big \vert  b^{(l)}(x) P_n^{(l,k-l)}(t,x) \Big \vert \\ & \leq 
%\left(\max_{l=1, \dots, k} {\Vert b^{(l)} \Vert}_{\infty}\right) \sum_{l=1}^{k} \sup_{(t,x) \in [0,T] \times (K \setminus A)} \Big \vert P_n^{(l,k-l)}(t,x) \Big \vert,
%\end{split}
%\end{equation*}
which also goes to zero for all $k \leq N$. Hence, ${\Vert (f - bP_n)^{(0,k)} \Vert}_{\infty} \lra 0$ as $\nto$ and the same holds for multi-indices $(k,0)$ with $k \leq N$ as $b$ is time independent. This concludes the proof. \qed
\epr

\bpr[of Theorem \ref{theo:weak:implies:mild}]
The proof is inspired by Theorem~9.15 in \cite{Peszat}. The key idea is that we can extend \eqref{weak:form:exact:Gauss} to test functions with a space \emph{and} a time variable. To be precise, for any $\psi_1, \psi_2 \in C_c^{\infty}([0,T] \times (-T,L+T))$, we will show that $\wt{\bbp}$-almost surely,
\begin{equation} \label{weak:form:exact:Gauss:ext}
\begin{split}
& \int_{\bbr} w(t,x) \psi_1(t,x) \, \dd x + \langle \theta_t, \psi_2(t, \cdot) \rangle \\ & \qquad \quad = 
\int_{0}^{t} \int_{\bbr} w(s,x) \left(\disfrac{\pd \psi_1}{\pd t}(s,x)  + 
\disfrac{\pd^2 \psi_2}{{\pd x}^2}(s,x)\right) \dd s \, \dd x + 
\int_{0}^{t} \left \langle \theta_s,  
\psi_1(s, \cdot) + \disfrac{\pd \psi_2}{\pd t}(s, \cdot)
\right \rangle \, \dd s  \\ & \qquad \qquad  \,  +
\int_{0}^{t} \int_{\bbr} \psi_2(s,x) f(w_{-}(s,x)) \, \wt{W}(\dd s, \dd x) \quad \textrm{for all} \quad t \leq T.
\end{split}
\end{equation}
First, we show \eqref{weak:form:exact:Gauss:ext} for special functions
\begin{equation} \label{spec:psi}
\Psi_i(t,x) = \varphi(t) \phi_i(x) \quad \textrm{with} \quad \varphi \in C^{\infty}([0,T]) \quad \textrm{and} \quad \phi_i \in C_c^{\infty}((-T,L+T)). 
\end{equation}
Using the integration by parts formula for càdlàg functions of Proposition~9.16 in \cite{Peszat} and taking into account the initial conditions of $w$ and $\theta$, we compute for all $t \leq T$,
\begin{equation} \label{NR10}
\begin{split}
& \int_{\bbr} w(t,x) \psi_1(t,x) \, \dd x + \langle \theta_t, \psi_2(t, \cdot) \rangle = \varphi(t) \left(\int_{\bbr} w(t,x) \phi_1(x) \, \dd x + 
\langle \theta_t, \phi_2 \rangle \right) \\ & \qquad = 
\int_{0}^{t} \varphi'(s) \left(\int_{\bbr} w(s,x) \phi_1(x) \, \dd x + 
\langle \theta_s, \phi_2 \rangle \right) \dd s + 
\int_{0}^{t} \varphi(s) \, \dd \left( \int_{\bbr} w(s,x) \phi_1(x) \, \dd x + 
\langle \theta_s, \phi_2 \rangle \right).
\end{split}
\end{equation}
Now the last integral process in \eqref{NR10} is indistinguishable from the process
\begin{equation} \label{NR20}
t \mapsto \int_{0}^{t} \varphi(s) \left(\int_{\bbr} w(s,x) \phi_2''(x) \, \dd x +  \langle \theta_s, \phi_1 \rangle \right) \dd s +
\int_{0}^{t} \int_{\bbr} \varphi(s) \phi_2(x) f(w_{-}(s,x)) \, \wt{W}(\dd s, \dd x),
\end{equation}
since its integrator equals the right-hand side of \eqref{weak:form:exact:Gauss} by assumption. Inserting \eqref{NR20} into \eqref{NR10} and recombining the functions $\psi_i$ as well as their derivatives yields exactly \eqref{weak:form:exact:Gauss:ext}. 

Next, we prove \eqref{weak:form:exact:Gauss:ext} for general $\psi_i \in C_c^{\infty}([0,T] \times (-T,L+T))$ by a density argument. Let $N_0 \in \bbn$ to be determined later in the proof. Using Lemma~\ref{lem:dense}, choose sequences $(\psi_i^n)_{n \in \bbn} \in C^{\infty}([0,T]) \otimes C_c^{\infty}((-T,L+T))$ such that $\psi_i^n$ converges to $\psi_i$ in $\sum_{\vert \al \vert \leq N_0} {\Vert \cdot \Vert}_{\infty, \al}$ with each
${\Vert \cdot \Vert}_{\infty, \al}$ as in \eqref{seminorm}. This implies uniform convergence in $[0,T]$ of each of the corresponding terms in \eqref{weak:form:exact:Gauss:ext} as we show in the following. (Note that by linearity, \eqref{weak:form:exact:Gauss:ext} readily holds for linear combinations of special functions \eqref{spec:psi}.)

Since $w$ is $\lpo$-càdlàg, and $\theta$ is càdlàg in $H_{-r}(\bbr)$, both processes are bounded and therefore,
\begin{equation*}
\sup_{t \leq T} \Bigg \vert \int_{\bbr} w(t,x) \psi_1(t,x) \, \dd x - \int_{\bbr} w(t,x) \psi_1^n(t,x) \, \dd x \Bigg \vert \leq 
C {\Vert w \Vert}_{\infty} {\Vert \psi_1 - \psi_1^n \Vert}_{\infty} \lra 0 \quad \textrm{as} \quad \nto
\end{equation*}
as well as
\begin{equation} \label{NR11}
\sup_{t \leq T} \Big \vert \langle \theta_t, \psi_2(t, \cdot) \rangle - \langle \theta_t, \psi_2^n(t, \cdot) \rangle \Big \vert \leq
\left(\sup_{t \leq T} {\Vert \theta_t \Vert}_{-r} \right)
\sup_{t \leq T} {\Vert  \psi_2(t, \cdot) -  \psi_2^n(t, \cdot) \Vert}_r < \infty.
\end{equation}
We now show that for any $r \geq 0$, $\psi_i^n \lra \psi_i$ in all ${\Vert \cdot \Vert}_{\infty, \al}$ with $\vert \al \vert \leq N_0$ and sufficiently large $N_0$ implies $\psi_i^n \lra \psi_i$ and $\pd_t \psi_2^n \lra \pd_t \psi_2$ in $\sup_{t \leq T} {\Vert \cdot \Vert}_r$ (thus convergence to 0 of all terms in \eqref{NR11}). For this, we use the well-known differential equation satisfied by the Hermite functions
\begin{equation} \label{Herm:EF}
h_q''(x) + (1 + 2q - x^2)h_q(x) = 0 \quad \textrm{for} \quad x \in \bbr \quad \textrm{and} \quad q \in \bbn.
\end{equation}
Let $q_0 \in \bbn$ be such that $\sqrt{1 + 2q} > L + T$ for all $q \geq q_0$. Then $1/\vert x^2 - (1 + 2q) \vert \leq 1/((1 + 2q) - (L + T)^2 )$ on $[-T, L + T]$ for all $q \geq q_0$. Let $\phi \in C_c^{\infty}((-T, L+T))$. Insert \eqref{Herm:EF} into $\langle \phi, h_q \rangle$ and use integration by parts twice, repeat $k$ times this procedure, apply then Hölder's inequality and the elementary inequality above to see that for all $q \geq q_0$ and $k \in \bbn$,
\begin{equation} \label{bound:skp:Herm}
\langle \phi, h_q \rangle^2 \leq
\frac{C }{{((1 + 2q) - (L + T)^2)}^{2k}} \left(\sum_{l=0}^{2k} {\Vert \phi^{(l)} \Vert}_{\infty}^2\right) 
\int_{-T}^{L + T} P_{4k}(x) \, \dd x
\end{equation}
with a polynomial $P_{4k}$ of degree $4k$ (the remaining details of these calculations are left to the reader). 
%Zwischenschritt:
%\begin{equation} \label{bound:skp:Herm}
%\begin{split}
%\langle \phi, h_q \rangle^2 & = \left(\int_{\bbr} \frac{\phi(x)}{x^2 - (1 + 2q)} h_q''(x) \, \dd x\right)^2 \\ & \leq
%\frac{C}{\left((1 + 2q) - (L + T)^2\right)^{2k}} 
%\left(\sum_{l=0}^{2k} {\Vert \phi^{(l)} \Vert}_{\infty}^2\right)
%\int_{-T}^{L + T} P_{2k}(x) \, \dd x,
%\end{split}
%\end{equation}
%\begin{equation*} 
%\begin{split}
%\int_{\bbr} \left( \frac{\phi''(x)}{x^2 - (1 + 2q)} - 
%\frac{4 x \phi'(x) + 2 \phi(x)}{(x^2 - (1 + 2q))^2} + \frac{8 x^2 \phi(x)}{(x^2 - (1 + 2q))^3} \right) h_q(x) \, \dd x 
%\end{split}
%\end{equation*}
Now choose $N_0$ such that $r - N_0 < - 3$ and infer from \eqref{bound:skp:Herm} that
%\begin{equation*}
%\begin{split}
%& {\Vert \psi_i^n(t,\cdot) - \psi_i(t,\cdot) \Vert}_r^2 \\ & \quad \leq 
%\sum_{q=0}^{q_0} (1 + 2q)^r {\langle \psi_i^n(t,\cdot) - \psi_i(t,\cdot), h_q \rangle}^2 +
%C \sum_{q=q_0+1}^{\infty} \frac{(1 + 2q)^r}{\left((1 + 2q) - (L + T)^2\right)^{2k_0}} 
%\sum_{l=0}^{2k_0} {\Vert \pd_x^{l} (\psi_i^n - \psi_i)  \Vert}_{\infty}^2
%\\ & \qquad \leq C
%\sum_{l=0}^{2k_0} {\Vert \pd_x^{l} (\psi_i^n - \psi_i)  \Vert}_{\infty}^2 
%\left( \sum_{q=0}^{q_0} (1 + 2q)^r + \sum_{q=q_0+1}^{\infty} \frac{(1 + 2q)^r}{\left((1 + 2q) - (L + T)^2\right)^{2k_0}}  \right)
%\end{split}
%\end{equation*}
\begin{equation} \label{NR16}
\begin{split}
& {\Vert \pd_t \psi_i^n(t,\cdot) - \pd_t \psi_i(t,\cdot) \Vert}_r^2 \\ & \qquad \leq C 
\left( \sum_{q=0}^{q_0} (1 + 2q)^r  +  \sum_{q=q_0+1}^{\infty} \frac{(1 + 2q)^r}{\left((1 + 2q) - (L + T)^2\right)^{N_0 - 2}}  \right) 
\sum_{l=0}^{N_0 - 2} {\Vert \pd_t \pd_x^{l} (\psi_i^n - \psi_i)  \Vert}_{\infty}^2
\end{split}
\end{equation}
for all $t \leq T$ and $n \in \bbn$. The series in \eqref{NR16} is finite since $r - N_0 < - 3$ and the last sum converges to 0 as $\nto$ by assumption on $(\psi_i^n)$, which proves the desired convergences. 

Finally, by Doob's inequality, It\={o}'s isometry, the Lipschitz continuity of $f$ and \eqref{ess:sup:2},
\begin{equation*}
\begin{split}
& \bbe \left[ \sup_{t \leq T} \Bigg \vert
\int_{0}^{t} \int_{\bbr} \psi_2(s,x) f(w_{-}(s,x)) \, \wt{W}(\dd s, \dd x) - \int_{0}^{t} \int_{\bbr} \psi_2^n(s,x) f(w_{-}(s,x)) \, \wt{W}(\dd s, \dd x) 
\Bigg \vert \right] \\ & \qquad \leq
\int_{0}^{T} \int_{\bbr} (\psi_2(s,x) - \psi_2^n(s,x))^2 \bbe \left[ f(w_{-}(s,x))^2 \right] \dd s \, \dd x \leq
C {\Vert \psi_2 - \psi_2^n \Vert}_{\infty}^2 \lra 0 \quad \textrm{as} \quad \nto.
\end{split}
\end{equation*}
Analogous arguments for the remaining terms of \eqref{weak:form:exact:Gauss:ext} finishes the density argument. 

We now choose two particular functions to be inserted in \eqref{weak:form:exact:Gauss:ext}. Fix $\phi \in C_c^{\infty}((0,L))$ as well as $t \leq T$, and define
\begin{equation} \label{spec:choi:psi}
\psi_1(s,y) = \frac{1}{2} \phi(y \pm (t-s)) \quad \textrm{and}  \quad
\psi_2(s,y) = \frac{1}{2} \int_{y - (t-s)}^{y + (t-s)} \phi(x) \, \dd x
\end{equation}
with $(s,y) \in [0,T] \times (-T, L + T)$.
%\begin{equation} \label{spec:choi:psi}
%\begin{split}
%\psi_1(s,y) & = \int_{\bbr} \phi(x) \frac{\dd  G}{\dd x}(t, \dd x; s, y) = \frac{1}{2} \phi(y \pm (t-s)) \,\,\, \textrm{and} \\
%\psi_2(s,y) & = \int_{\bbr} \phi(x) G(t,x;s,y) \, \dd x = \frac{1}{2} \int_{y - (t-s)}^{y + (t-s)} \phi(x) \, \dd x, \,\,\, (s,y) \in [0,T] \times (-T, L + T).
%\end{split}
%\end{equation}
Then $\psi_1, \psi_2 \in C_c^{\infty}([0,T] \times (-T,L+T))$. In addition, $\psi_1(t,y) = \phi(y)$ and $\psi_2(t,y) = 0$ for all $y \in \bbr$, and straightforward calculus yields
\begin{equation} \label{deriv:spec:psi}
\psi_2^{(1,0)}(s,y) = - \psi_1(s,y) \quad \textrm{and} \quad \psi_2^{(2,0)}(s,y) = \psi_2^{(0,2)}(s,y).
\end{equation}
% \disfrac{\pd \Psi_2}{\pd s}(s,y) = - \Psi_1(s,y) \quad \textrm{and} \quad \disfrac{\pd^2 \Psi_2}{\pd s^2}(s,y) = \disfrac{\pd^2 \Psi_2}{\pd y^2}(s,y).
The freedom we have to choose two different functions in \eqref{spec:choi:psi} is another reason why we considered the weak formulation \eqref{weak:form:formal} of the stochastic wave equation in this work: By \eqref{deriv:spec:psi}, the first two integrals on the right-hand side of \eqref{weak:form:exact:Gauss:ext} \emph{vanish}, and \eqref{weak:form:exact:Gauss:ext} yields at time point $t$
%Insert $\psi_1$ and $\psi_2$ into \eqref{weak:form:exact:Gauss:ext} to see that the first two integrals on the right-hand side \emph{vanish} and that \eqref{weak:form:exact:Gauss:ext} yields at time point $t$:
\begin{equation*}
\int_{\bbr} w(t,x) \phi(x) \, \dd x = \int_{0}^{t} \int_{\bbr} \psi_2(s,y) f(w_{-}(s,y)) \, \wt{W}(\dd s, \dd y) \quad \wt{\bbp}\textrm{-almost surely},
\end{equation*}
which, recalling \eqref{dist:deriv:G} and using the stochastic Fubini theorem, has the equivalent form 
\begin{equation} \label{NR12}
\int_{\bbr} \left( w(t,x) - \int_{0}^{t} \int_{\bbr} G(t,x;s,y) f(w_{-}(s,y)) \, \wt{W}(\dd s, \dd y) \right) \phi(x) \, \dd x = 0 \quad \wt{\bbp} \textrm{-almost surely},
\end{equation}
and this holds for all $\phi \in C_c^{\infty}((0,L))$ and $t \leq T$. We can now infer the first claim of the theorem. Denote by $Z_t(x)$ the random field in parenthesis in \eqref{NR12} with $(t,x) \in [0,T] \times [0,L]$. It is easy to see that \eqref{ess:sup:2} implies $Z_t \in L^2([0,L])$ for all $t \leq T$.
%This follows from the assumptions on $w$, It\={o}'s isometry and the Lipschitz continuity of $f$:
%\begin{equation*}
%\begin{split}
%\bbe \left[ \int_{0}^{L} Z_t^2(x) \, \dd x \right] & \leq 
%2 \int_{0}^{L} \left(\bbe \left[ w(t,x)^2 \right] + \bbe \left[ \int_{0}^{t} \int_{\bbr} G^2(t,x;s,y) f^2(w_{-}(s,y)) \, \dd s \, \dd y \right] \right) \dd x \\ & \leq
%C \left(\esssup_{(t,x) \in [0,T] \times [-T,L+T]} \bbe \left[ \big \vert w(t,x) \big \vert^2 \right] + 1\right) < \infty.
%\end{split}
%\end{equation*}
For any $\epsilon > 0$ and $t \leq T$, consider the mollified random field $J_{\epsilon}(Z_t)$ on $[0,L]$ defined exactly as in (1.8) of Chapter 10 in \cite{Friedman}. By Lemma~3 of that chapter, $J_{\epsilon}(Z_t) \lra Z_t$ in $L^2([0,L])$ as $\epsilon \to 0$. Consequently, we can choose a sequence $(\epsilon_l)_{l \in \bbn}$ converging to 0 such that $\wt{\om}$-wise,
\begin{equation} \label{alm:sure:conv:moll}
J_{\epsilon_l}(Z_t) \lra Z_t \quad \Leb_{[0,L]} \textrm{-almost everywhere} \quad \textrm{as} \quad l \to \infty.
\end{equation}
Now for any fixed $y \in (0,L)$, the support of the function 
$\rho((y - \cdot)/\epsilon_l)/ \epsilon_l^2 $ used to mollify $Z_t$ will be contained in $(0,L)$ if $\epsilon_l$ is sufficiently small and, hence, \eqref{NR12} applies to $J_{\epsilon_l}(Z_t)(y)$ for those $\epsilon_l$ with $\phi$ being chosen as $\rho((y - \cdot)/\epsilon_l)/ \epsilon_l^2 $. Combining this with \eqref{alm:sure:conv:moll} has the following outcome: For all $t \leq T$ and almost all $y \in (0,L)$, $Z_t(y) = 0$ $\wt{\bbp}$-almost surely. We deduce that $w_{-}$ satisfies the mild formulation of \eqref{WAVE:Gauss:2} almost everywhere on $[0,T] \times [0,L]$.

Let $\wt{u}$ be the continuous mild solution to \eqref{WAVE:Gauss:2} on $(\wt{\Om}, \wt{\calf}, \wt{\bfil}, \wt{\bbp})$. We have 
\begin{equation*}
\bbe \left[ \Big \vert \wt{u}(t,x) - w_{-}(t,x) \Big \vert^2 \right] \leq 
\int_{0}^{T} \int_{\bbr} G^2(t,x;s,y) \bbe \left[ \Big \vert \wt{u}(s,y) - w_{-}(s,y) \Big \vert^2 \right] \dd s \, \dd y \quad \textrm{a.e.}
\end{equation*}
and by Lemma 6.4 (3) in \cite{CC2}, $\bbe \left[ \vert \wt{u}(t,x) - w_{-}(t,x) \vert^2 \right] = 0$ for almost all $(t,x) \in [0,T] \times [0,L]$. It follows that $\wt{\bbp}$-almost surely, $w$ and $\wt{u}$ agree almost everywhere on $[0,T] \times [0,L]$ and therefore, since $w$ is $\lpo$-càdlàg, they are indistinguishable and $w$ is actually continuous on $[0,T] \times [0,L]$. Finally, by the usual computations, we obtain for all $(t,x) \in [0,T] \times [0,L]$,
\begin{equation*}
\begin{split}
w(t,x) = \wt{u}(t,x) & = \int_{0}^{t} \int_{\bbr} G_{t-s}(x,y) f(\wt{u}(s,y)) \, \wt{W}(\dd s, \dd y) = \int_{0}^{t} \int_{\bbr} G_{t-s}(x,y) f(w(s,y)) \, \wt{W}(\dd s, \dd y)
\end{split}
\end{equation*}
$\wt{\bbp}$-almost surely.

For the second claim of the theorem, we define new functions $\psi_1, \psi_2$ by
\begin{equation*}
\psi_1(s,y) = \frac{1}{2} \left( \phi'(y+(t-s)) - \phi'(y - (t-s)) \right) \quad \textrm{and} \quad
\psi_2(s,y) = \frac{1}{2} \phi(y \pm (t-s))
\end{equation*}
with fixed $\phi \in C_c^{\infty}((0,L))$ and $t \leq T$, for all $(s,y) \in [0,T] \times (-T, L + T)$. Again we have $\psi_1, \psi_2 \in C_c^{\infty}([0,T] \times (-T, L + T))$. Since $w_{-} = w$, by straightforward calculus and again with \eqref{dist:deriv:G}, inserting $\psi_1, \psi_2$ into \eqref{weak:form:exact:Gauss:ext} yields at time point $t$ exactly \eqref{NR18}. This concludes the proof. \qed
\epr

\subsection*{Acknowledgements}
TD cordially thanks Carsten Chong as well as Claudia Klüppelberg for inspiring discussions and valuable advice. TD's research is partially supported by the Deutsche Forschungsgemeinschaft, project number KL 1041/7-1.

\addcontentsline{toc}{section}{References}
\bibliographystyle{plainnat}
\bibliography{bib-Stochastic_Wave_Equation}

\end{document}